%% file: BlaSU21b.tex
\documentclass[11pt,reqno]{amsart}

\usepackage[english]{babel}
\usepackage[utf8]{inputenc}
\usepackage[T1]{fontenc}
\usepackage{lmodern} \normalfont 
\DeclareFontShape{T1}{lmr}{bx}{sc} { <-> ssub * cmr/bx/sc }{}
\usepackage{microtype}	
\usepackage{xspace}

\usepackage{amssymb}
\usepackage{amsmath}
\usepackage{amsthm}
\usepackage{mathtools}
\usepackage{mathrsfs}
\usepackage{accents} 
\usepackage{etoolbox}
\usepackage{siunitx}
\usepackage{dsfont} 
\usepackage{graphicx}
\usepackage{color}
\usepackage[dvipsnames]{xcolor}
\usepackage{tikz}
\usetikzlibrary{calc,positioning,shapes}
\usetikzlibrary{patterns,decorations.pathmorphing,decorations.markings}
\usepackage{pgfplots}
\pgfplotsset{compat=1.17}

\usepgfplotslibrary{groupplots}

\pgfplotsset{
        colormap={parula}{
            rgb255=(53,42,135)
            rgb255=(15,92,221)
            rgb255=(18,125,216)
            rgb255=(7,156,207)
            rgb255=(21,177,180)
            rgb255=(89,189,140)
            rgb255=(165,190,107)
            rgb255=(225,185,82)
            rgb255=(252,206,46)
            rgb255=(249,251,14)
        },
}

\usepackage[margin=10pt,font=small,labelfont=bf,labelsep=endash]{caption}
\usepackage{subcaption}

\usepackage{booktabs}
\usepackage{paralist}
\usepackage{soul}

\numberwithin{equation}{section}

\usepackage{algorithm}
\usepackage{algorithmicx}
\usepackage{algpseudocode}

\usepackage{cite}
\usepackage{multirow} 
\usepackage{enumitem} 
\setlist[enumerate]{label=(\roman*)}

\usepackage[colorlinks,citecolor=black,urlcolor=black,linkcolor=black]{hyperref}
\usepackage[nameinlink,noabbrev]{cleveref}

\theoremstyle{plain}
\newtheorem{theorem}{Theorem}[section]

\newtheorem{remark}[theorem]{Remark}

\newtheorem{example}[theorem]{Example}

\input{def.tex}
\newcommand{\lineWidth}{1.2}
\newcommand{\markerSize}{3}

\newcommand{\beforeSwitchTag}{(i)\xspace}
\newcommand{\beforeSPODTag}{(ii--a)\xspace}
\newcommand{\sPODTag}{(ii--b)\xspace}

\title[Efficient Wildland Fire Simulation via Nonlinear MOR]{Efficient Wildland Fire Simulation via Nonlinear Model Order Reduction}
\author[F.~Black, P.~Schulze, and B.~Unger]{Felix Black ${}^\dagger$ \and Philipp Schulze${}^\dagger$ \and Benjamin Unger${}^\star$}
\address{${}^{\dagger}$  Institute of Mathematics MA\,{}4-4, Technical University Berlin, Stra\ss e des 17.~Juni 136, 10623 Berlin, Germany}
\email{\{black,pschulze\}@math.tu-berlin.de}
\address{${}^{\star}$ Stuttgart Center for Simulation Science (SC SimTech), University of Stuttgart, Universit\"{a}tsstr.~32, 70569 Stuttgart, Germany}
\email{benjamin.unger@simtech.uni-stuttgart.de}

\date{\today}

\begin{document}

\begin{abstract}
	We propose a new hyper-reduction method for a recently introduced nonlinear model reduction framework based on dynamically transformed basis functions and especially well-suited for advection-dominated systems.
Furthermore, we discuss applying this new method to a wildland fire model whose dynamics feature traveling combustion waves and local ignition and is thus challenging for classical model reduction schemes based on linear subspaces.
The new hyper-reduction framework allows us to construct parameter-dependent reduced-order models (ROMs) with efficient offline/online decomposition.
The numerical experiments demonstrate that the ROMs obtained by the novel method outperform those obtained by a classical approach using the proper orthogonal decomposition and the discrete empirical interpolation method in terms of run time and accuracy.
\end{abstract}

\maketitle

%
%

\section{Introduction}

\begin{figure}
	\centering
	\begin{subfigure}[t]{.32\linewidth}
		\begin{tikzpicture}
			\pgfplotsset{width=6cm,height=6cm}
			\begin{groupplot}[group style={group size=1 by 1, horizontal sep=2.5cm}]
				\nextgroupplot[enlargelimits=false,axis on top=false,xlabel=$x$,ylabel=$t$,xtick=\empty,ytick=\empty,ylabel style={rotate=-90}]
				\addplot graphics [xmin=0,xmax=1,ymin=0,ymax=1]{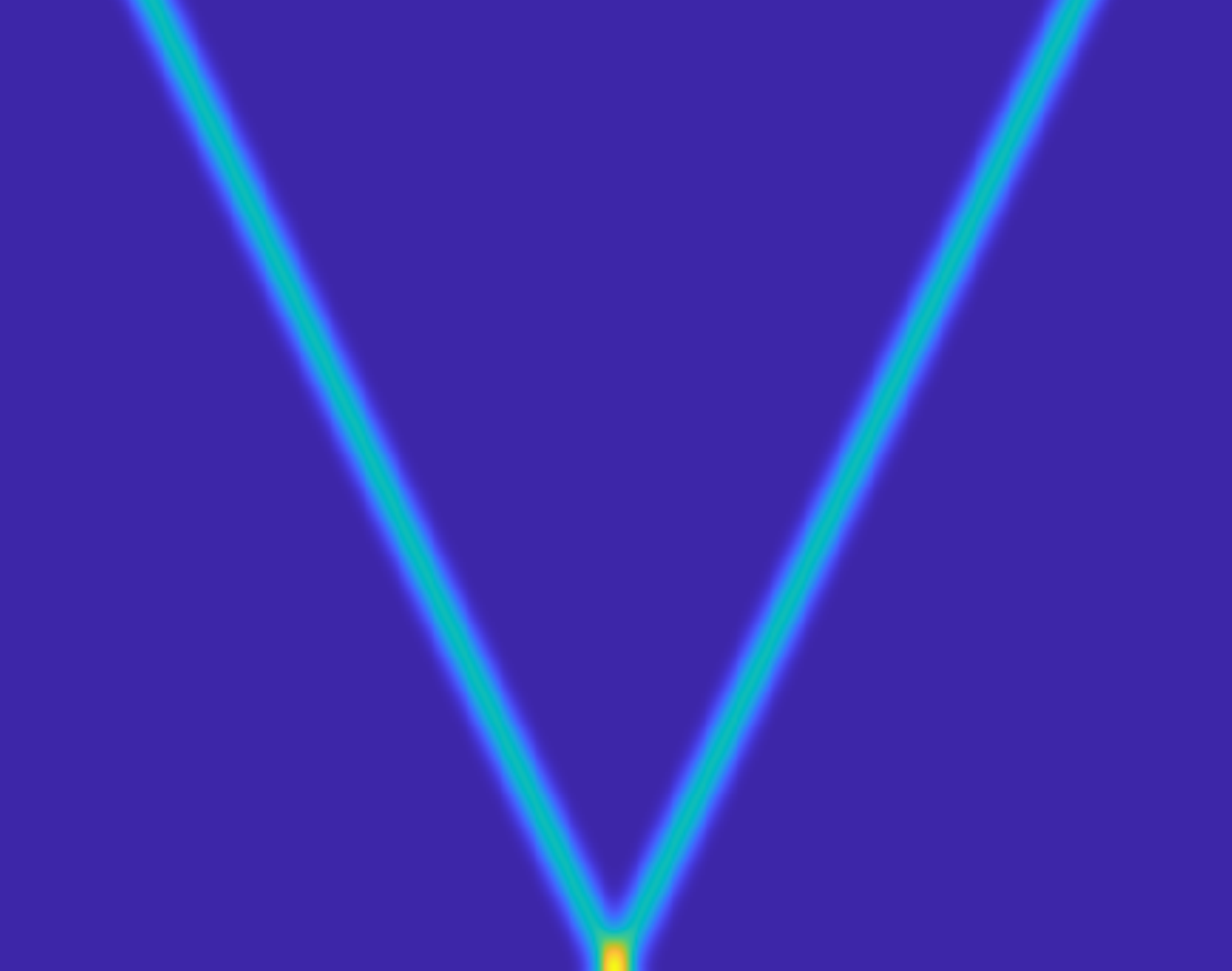};
			\end{groupplot}
		\end{tikzpicture}
		\caption{Wave equation}
		\label{fig:snapshot:waveEquation}
	\end{subfigure}\hfill
	\begin{subfigure}[t]{.32\linewidth}
		\begin{tikzpicture}
			\pgfplotsset{width=6cm,height=6cm}
			\begin{groupplot}[group style={group size=1 by 1, horizontal sep=2.5cm}]
				\nextgroupplot[enlargelimits=false,axis on top=false,xlabel=$x$,xtick=\empty,ytick=\empty,ylabel style={rotate=-90},ylabel=$t$]
				\addplot graphics [xmin=0,xmax=1,ymin=0,ymax=1]{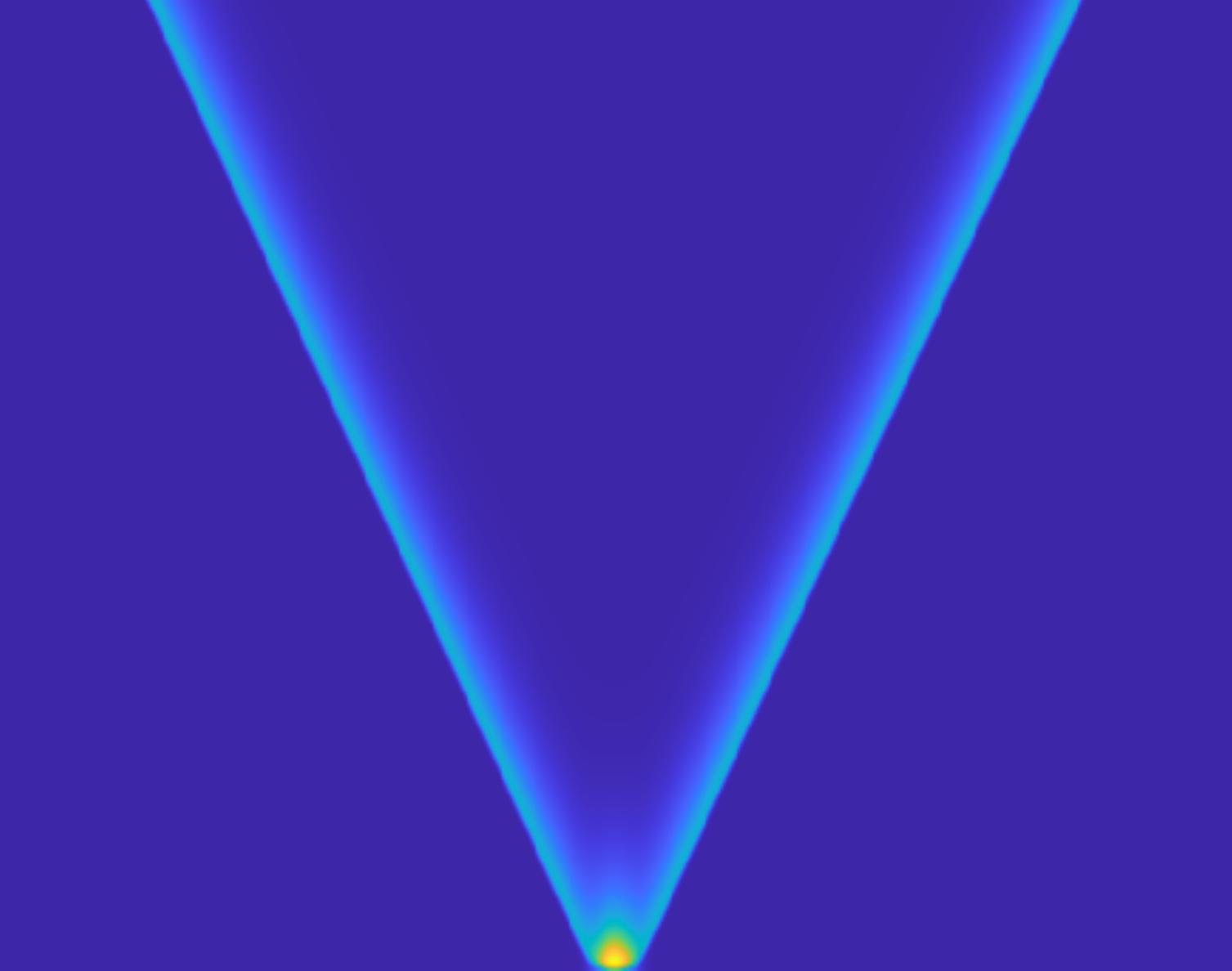};
			\end{groupplot}
		\end{tikzpicture}
		\caption{Temperature}	
		\label{fig:snapshot:wildlandFire:temperature}
	\end{subfigure}\hfill
	\begin{subfigure}[t]{.32\linewidth}
		\begin{tikzpicture}
			\pgfplotsset{width=6cm,height=6cm}
			\begin{groupplot}[group style={group size=1 by 1, horizontal sep=2.5cm}]
				\nextgroupplot[enlargelimits=false,axis on top=false,xlabel=$x$,xtick=\empty,ytick=\empty,ylabel style={rotate=-90},ylabel=$t$]
				\addplot graphics [xmin=0,xmax=1,ymin=0,ymax=1]{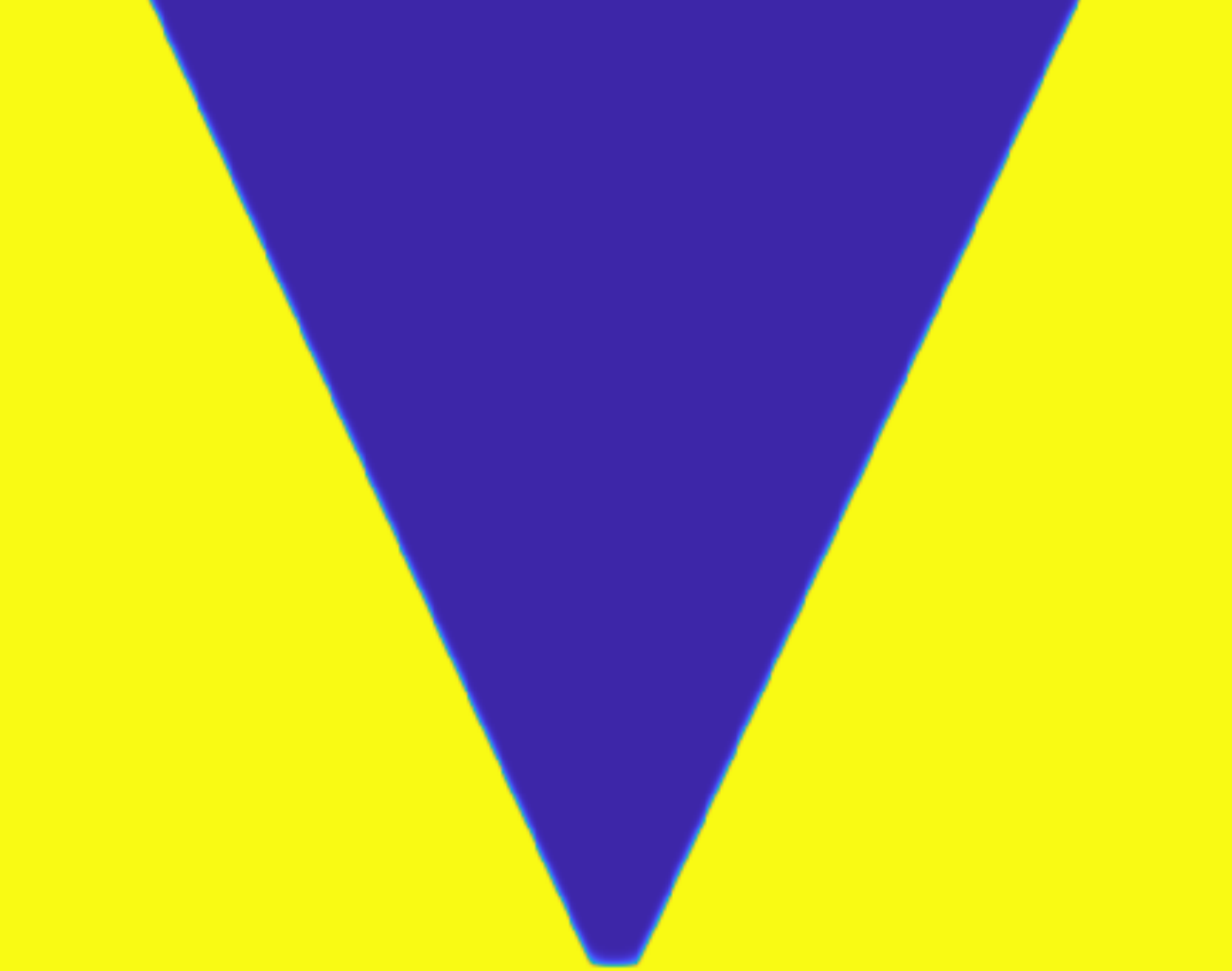};
			\end{groupplot}
		\end{tikzpicture}
		\caption{Supply mass fraction}	
		\label{fig:snapshot:wildlandFire:smf}
	\end{subfigure}
	\caption{Comparison of (dimensionless) snapshot matrices of a particular solution of the wave equation and a
	solution (temperature and supply mass fraction) of a wildland fire simulation}
	\label{fig:snapshot}
\end{figure}

Accurate modeling of the behavior of wildland fires is an important but challenging task. The models typically involve different physical domains,  see \cite{Per98,ManBBCDKV08} for an overview, or in parts even data-integrated models~\cite{LatHL20}.  
Besides the multi-physics character of the problem, different scales of the involved quantities further complicate the simulation of the coupled models.  
Consequently, high-fidelity models are computationally costly, conspiring against quick forecasting results required in emergency situations. To reduce the computational time while retaining the complexity of the model, model order reduction (MOR) is a promising approach that has been applied successfully to many different applications. For an overview of MOR methods and applications, we refer to \cite{Ant05,BenCOW17,HesRS16,QuaR14,QuaMN16,AntBG20}.  
Standard projection-based MOR methods such as proper orthogonal decomposition (POD), balanced truncation, or rational interpolation, perform MOR by choosing suitable low-dimensional subspaces of the solution space and then project the system dynamics onto these subspaces using a Petrov--Galerkin approach. 
The best subspace of a given dimension is classified by the Kolmogorov $n$-widths~\cite{Kol36,Pin85} (respectively the Hankel singular values \cite{UngG19}), thus providing a lower bound on the approximation error. 
Although in many applications the Kolmogorov $n$-widths decay exponentially~\cite{MadPT02b}, thus enabling MOR to succeed, this is not true for every application. 
A prominent example is the linear wave equation \cite{GreU19}, whose (dimensionless) snapshot matrix based on a Gaussian initial condition is presented in \Cref{fig:snapshot:waveEquation}. Comparing this snapshot matrix with the solution of a wildland fire model that couples the temperature (illustrated in \Cref{fig:snapshot:wildlandFire:temperature}) with the available supply mass fraction (cf.~\Cref{fig:snapshot:wildlandFire:smf}) suggests that the $n$-widths do not decay exponentially for this model either.
In particular, we cannot expect standard MOR methods to produce an accurate and, at the same time, low-dimensional approximation for wildland fire simulations. To overcome slowly decaying Kolmogorov $n$-widths, MOR methods for transport-dominated phenomena have emerged over the last couple of years.  Let us refer to \cite{OhlR13,TadPQ15,RimML18,LeeC19,CagMS19,NonBRM19,BlaSU20,Peh20,Tad20,Wel20,KraSR21} and the references therein to name just a few. 
As a general strategy, many of these methods use some kind of nonlinear projection framework to allow the ansatz space to adapt along with the solution. 

In this paper, we apply the method presented in \cite{BlaSU20}, which can be understood as a generalization of the
moving finite element method \cite{MilM81} to general ansatz functions, to a simple wildland fire model introduced in
\cite{ManBBCDKV08}. 
This showcases the potential of MOR methods tailored to transport phenomena in real-time wildland fire simulations.  
Our framework builds on an efficient (offline) decomposition \cite{ReiSSM18,Rei21,SchRM19} of already computed or measured solution data and constructs the reduced-order model (ROM) via residual minimization as detailed in~\cite{LeeC19,BlaSU20}. 
Due to the nonlinearity of the projection method, the resulting ROM still depends on the original dimension, which conspires against an efficient simulation. Besides, the wildland fire model itself is nonlinear due to a (modified) Arrhenius law, which provides an additional challenge for MOR, even if standard methods are applied. To avoid scaling with the full dimension, it is common practice to work with a further approximation of the nonlinearity, typically referred to as hyper-reduction.  A prominent hyper-reduction methods are the empirical interpolation method (EIM) \cite{BarMNP04} or its discrete analog (DEIM) \cite{ChaS10}.  In this paper, we extend this approach to cope with dynamically transformed modes. Our main results are the following:
\begin{itemize}
	\item In \cref{subsec:shiftedDEIM}, we extend DEIM by approximating the nonlinearity of the full-order model (FOM) by a linear combination of dynamically transformed ansatz functions to account for the advective transport within the system.
		Furthermore, we discuss how the state-dependent ROM coefficient matrices can be replaced with cheap to evaluate approximants, see \cref{sec:innerProducts}.
		Altogether, the proposed methodology allows for constructing parameter-dependent ROMs while achieving an efficient offline/online decomposition.
	\item Based on the new approach, we construct low-dimensional parameter-dependent ROMs for a wildland fire model, see \cref{sec:numerics}. 
	Depending on the initial condition, the model inherits complex dynamics that do not allow for a simple separation of transported and non-transported effects, thus rendering this a challenging benchmark problem. 
	Following ideas from \cite{DihDH11}, we propose a switching strategy, using our method solely in the transport-dominated regime, see \cref{sec:Gaussian} for further details. 
		The constructed ROMs allow for accurate predictions within the parameter space and prove to be faster and more accurate than POD-DEIM ROMs, which are based on classical linear subspace approximations.
\end{itemize}

To render the manuscript self-contained, we first present the wildfire model from \cite{ManBBCDKV08} in \cref{sec:wildfireModel}. 
Classical projection-based MOR for a general nonlinear dynamical system is reviewed in \cref{subsec:standardMOR}. Hyper-reduction via DEIM is briefly summarized in \cref{subsec:onlineOfflinePOD}.
In \cref{subsec:shiftedPOD}, we recall the shifted POD (sPOD)~\cite{ReiSSM18} as a powerful offline decomposition method together with the construction of a corresponding ROM based on dynamically transformed modes \cite{BlaSU20}. 

\subsection*{Notation}

The space of real $m\times n$ matrices is denoted by $\R^{m\times n}$, and for the transpose of a matrix $A$, we write $A^\top$.
Furthermore, we use the symbol $I_n$ for the identity matrix of size $n$ and denote its entries by the Kronecker delta $\delta_{ij} = [I_n]_{ij}$ for $i,j=1,\ldots, n$.
For diagonal matrices and blockdiagonal matrices, we use the abbreviations
\begin{equation*}
	\diag(a_1,\ldots,a_n) \vcentcolon=
	\begin{bmatrix}
		a_1 	& 				& \\
				& \ddots 	& \\
				&				& a_n
	\end{bmatrix}
	,\quad \blkdiag(A_1,\ldots, A_n) \vcentcolon=
	\begin{bmatrix}
		A_1 	& 				& \\
				& \ddots 	& \\
				&				& A_n
	\end{bmatrix}
	,
\end{equation*}
respectively, where $a_1,\ldots,a_n$ denote scalars and $A_1,\ldots, A_n$ matrices of arbitrary size.
In addition, the Kronecker product between two matrices $A$ and $B$ is denoted by $A\otimes B$.
The symbols $\langle \cdot, \cdot\rangle$ and~$\lVert \cdot \rVert$ are used for the Euclidean inner product and norm, respectively.
Furthermore, the space of square-integrable functions over the interval $(a,b)\subset\R$ is denoted by $\Ltwo(a,b)$. 

\section{Problem Setting -- Wildfire Model}
\label{sec:wildfireModel}

We use the model from \cite{ManBBCDKV08}. 
Considering a Lipschitz domain $\Omega\in\mathbb{R}^d$ with $d\in\{1,2\}$, a finite time horizon $0<\finalTime<\infty$, and setting $\mathbb{T} \vcentcolon= [0,\finalTime]$, one wants to determine the temperature $T\colon \mathbb{T} \times \Omega \to \mathbb{R}$ and fuel supply mass fraction $S\colon \mathbb{T}\times \Omega\to\mathbb{R}$ satisfying the coupled partial differential equation (PDE)
\begin{subequations}
	\label{eqn:wildlandfire:strong}
	\begin{align}
		\label{eqn:energyBalance}
		\partial_t T &= \nabla\cdot\,(k \nabla T) - v\cdot\nabla T + \tempRise\big(Sr(T,\arrhenius,T_{\mathrm{a}}) - \heatTransfer(T-T_{\mathrm{a}})\big),\\
		\label{eqn:partialMassBalance}
		\partial_t S &= -\fuelRate Sr(T,\arrhenius,T_{\mathrm{a}}),
	\end{align}
\end{subequations}
with reaction rate constant $\fuelRate r(T,\arrhenius,T_{\mathrm{a}})$ given by a modified Arrhenius law of the form
\begin{equation}
	\label{eqn:reactionRate}
	\fuelRate r(T,\arrhenius,T_{\mathrm{a}}) \vcentcolon= \begin{cases}
		\fuelRate\exp\big(-\tfrac{\arrhenius}{T-T_{\mathrm{a}}}\big), & \text{if } T > T_{\mathrm{a}},\\
		0, & \text{otherwise.}
	\end{cases}
\end{equation}
Here, $k$ denotes the thermal diffusivity, $\tempRise$ a proportionality constant representing the ratio of temperature rise per second and reaction rate, $\arrhenius$ the proportionality coefficient in the modified Arrhenius law, $\heatTransfer$ a scaled heat transfer coefficient for the heat loss to the environment, $\fuelRate$ the pre-exponential factor, $T_{\mathrm{a}}$ the ambient temperature, and $v$ the wind speed (given by atmospheric data or model). 
An overview of the variables, including the units, is provided in \Cref{tab:wildandfireQuantities}. For notational convenience, we perform a variable substitution and work with the relative temperature $T-T_{\mathrm{a}}$ instead of the temperature, which formally results in $T_{\mathrm{a}}$.

\begin{table}
	\centering
	\caption{Variables in the wildland fire model}
	\label{tab:wildandfireQuantities}
	\resizebox{0.95\textwidth}{!}{  
	\begin{tabular}{lrl@{\hspace{2em}}lrl}
		\toprule
		name & symbol & unit & name & symbol & unit\\\midrule
		temperature & $T$ & \si{\kelvin} &
		supply mass fraction & $S$ &\\
		thermal diffusivity & $k$ & \si{\square\metre\per\second} &
		temperature rise per second & $\tempRise$ & \si{\kelvin\per\second}\\
		proportionality coefficient & $\arrhenius$ & \si{\kelvin} &
		scaled heat transfer coefficient & $\heatTransfer$ & \si{\per\kelvin} \\
		pre-exponential factor & $\fuelRate$ & \si{\per\second} &
		ambient temperature & $T_{\mathrm{a}}$ & \si{\kelvin}\\
		wind speed & $v$ & \si{\metre\per\second}\\
		\bottomrule
	\end{tabular}}
\end{table}

For the derivation of the model, fire in a ground layer with finite, small thickness is considered. The fire layer consists of fuel and air above the fuel. 
The (initially two-dimensional) model is then derived from physical principles, namely from conservation of energy in the fire layer and a balance equation for the fuel supply. For a detailed discussion of the derivation, we refer to \cite[Section 4]{ManBBCDKV08}.
One of the modeling assumptions is that the only heat source stems from the heat generated by the chemical reaction of burning.
Furthermore, it is assumed that heat loss to the atmosphere is only caused by convection and radiation. In contrast, the short-range heat transfer due to radiation and turbulence is modeled via diffusion. 
More precisely, the individual terms in the model correspond to the following physical effects: 
\begin{itemize}
	\item The term $\nabla\cdot \left( k \nabla T \right)$ is a diffusion term modeling the short-range heat transfer by radiation and turbulence. 
	\item The term $v\cdot\nabla T$ accounts for the advective heat transfer caused by the wind with wind speed $v$.
	\item The loss of fuel due to burning is modeled by the product of the supply mass fraction and the reaction rate constant, which is modeled with the modified Arrhenius law~\eqref{eqn:reactionRate}.
	\item	 The term $\tempRise\heatTransfer \left( T - T_{\mathrm{a}} \right)$ models heat lost to the atmosphere due to convection. 
		It is assumed that heat loss to the atmosphere due to radiation is negligible compared to heat loss due to convection, such that the previous term is sufficient to account for both effects. 
\end{itemize}
We remark that the reason for the usage of the modified Arrhenius law lies in the so-called \emph{cold boundary difficulty} \cite{BerLR91}: If a combustion model does not enforce the absence of combustion at ambient temperature, the fuel supply at any given point in space is slowly depleted due to this permanent combustion effect, even before the combustion wave started by some local ignition process reaches this point. 
Note that the reaction rate constant with modified Arrhenius law is smooth in $T$. 

Physical effects that are not accounted for in this model are the effect of fuel disappearance on the total heat capacity of the fuel layer, storage of heat in the ground, chemical kinetics of intermediate products of combustion, the fine scale dynamics of the fire on the surface of and within the wood (fuel), and evaporation of moisture, cf.~\cite{ManBCK09}.
The model also does not include any two-way coupling with an atmospheric model. 
Lastly, note that the Arrhenius law is in general only valid for premixed flames~\cite{FerSLW06} with a sufficient oxygen supply and ignores fuel surface effects.
 
\section{Preliminaries}

\subsection{Projection-based MOR and proper orthogonal decomposition}
\label{subsec:standardMOR}
To simplify the following discussion, we assume that we have a finite-dimensional approximation of the coupled PDE~\eqref{eqn:wildlandfire:strong} available, obtained via semi-discretization in space, given by
\begin{equation}
	\label{eqn:truthModel}
	\begin{aligned}
		\dot{T}_h &= A_T(k,v,\tempRise,\heatTransfer)T_h + \tempRise f(S_h,T_h,\arrhenius),\\
		\dot{S}_h &= -\fuelRate f(S_h,T_h,\arrhenius),
	\end{aligned}
\end{equation}
with matrix $A_T$ and nonlinear function $f$. 
We assume that the approximation error is sufficiently small compared to the model reduction error to come, and refer to~\eqref{eqn:truthModel} as the \emph{truth model}. Combining the discretized temperature $T_h$ and supply mass fraction $S_h$ in a single state variable $\state$, and the remaining variables in the parameter vector $\param$, we are thus faced with the task of determining $\state\colon\timeInt\times\paramSet\to\R^{\stateDim}$ for some $\param\in\paramSet$ with compact parameter set $\paramSet$ satisfying
\begin{equation}
	\label{eqn:FOM}
	\left\{\begin{aligned}\quad
		\dot{\state}(t;\param) &= A(\param)\state(t;\param) + F(\state(t;\param),\param), & t\in\timeInt,\\
		\state(0;\param) &= \state^0.
	\end{aligned}\right.
\end{equation}
Depending on the spatial discretization method, the models in~\eqref{eqn:truthModel} and~\eqref{eqn:FOM} may have symmetric positive definite mass matrices in front of the time-derivatives. 
To simplify notation, we decided to ignore the mass matrices and emphasize that the discussion to follow can easily be extended to the case where the mass matrix is not equal to the identity matrix, by using a suitable weighted inner product.
As a further simplification and in agreement with the wildland fire model~\eqref{eqn:wildlandfire:strong}, we assume that the matrix~$A(\param)$ has an affine decomposition with respect to the parameter, i.e., for some $k\in\N$ and $\nu=1,\ldots,k$ there exist real functions $q_\nu\colon\paramSet\to\mathbb{R}$ and matrices $A_\nu\in\R^{\stateDim\times\stateDim}$ satisfying
\begin{equation}
	\label{eqn:affineDecomposition}
	A(\param) = \sum_{\nu=1}^k q_\nu(\param)A_\nu.
\end{equation}
The task of projection-based MOR is to identify a suitable low-dimensional subspace $\MORspace\subseteq \R^{\stateDim}$, parametrized via an orthonormal basis $(\varphi_1,\ldots,\varphi_\ROMdim)\in \MORspace^\ROMdim$, such that the solution of~\eqref{eqn:FOM} approximately evolves within $\MORspace$ for any $\param\in\paramSet$. 
Assembling the matrix $V = \begin{bmatrix}
	\varphi_1 & \ldots & \varphi_\ROMdim\end{bmatrix}\in\R^{\stateDim\times\ROMdim}$ and performing a Galerkin projection, we define
\begin{align*}
	\Ared_\nu &\vcentcolon= V^\top A_\nu V, & 
	\Ared(\param) &\vcentcolon= \sum_{\nu=1}^k g_\nu(\param)\Ared_\nu, & 
	\fred(\ROMstate,\param) &\vcentcolon= V^\top F(V\ROMstate,\param)\in\R^{\ROMdim},
\end{align*}
yielding the ROM
\begin{equation}
	\label{eqn:ROM:2}
	\left\{\begin{aligned}\quad
		\dot{\ROMstate}(t;\param) &= \Ared(\param)\ROMstate(t;\param) + \fred(\ROMstate(t;\param),\param),&
		t\in\timeInt, \\
		\ROMstate(0;\param) &= \ROMstate^0,
	\end{aligned}\right.
\end{equation}
with projected initial value $\ROMstate^0\vcentcolon= V^T\state^0$. The approximation obtained by the ROM is then given as 
\begin{equation}
	\label{eqn:MORansatz}
	\state(t;\param) \approx \ROMapprox(t;\param) \vcentcolon= V\ROMstate(t;\param).
\end{equation}
The question that remains to be answered is how to determine the approximation space~$\MORspace$, respectively the matrix $V$.  
If solution trajectories $\state(\cdot;\param_\sigma)$ are available for $\sigma=1,\ldots,\ell$ with parameter values $\param_\sigma\in\paramSet$, then the \emph{proper orthogonal decomposition} (POD), see for instance \cite{GubV17}, constructs the reduced basis by solving the minimization problem
\begin{equation}
	\label{eqn:PODminimization}
	\left\{\begin{aligned}
		&\min \frac{1}{2} \sum_{\sigma=1}^\ell \int_0^{\finalTime} \bigg\| \state(t;\param_\sigma) - \sum_{i=1}^\ROMdim \langle\state(t;\param_\sigma),\varphi_i\rangle \varphi_i \bigg\|^2 \mathrm{d}t\\
		&\mathrm{s.t.}~\varphi_i\in\R^{\stateDim} \quad\text{and}\quad \langle\varphi_i,\varphi_j\rangle =
		\delta_{ij}\quad\text{for } i,j=1,\ldots,\ROMdim.
	\end{aligned}\right.
\end{equation}
The POD basis can be computed efficiently via a truncated \emph{singular value decomposition} (SVD) of the snapshot matrix, given by concatenating the solution vectors at different time instances, which are typically given from a time discretization. 
We emphasize that the approximation quality throughout the parameter set $\paramSet$ will depend on the choice of the \emph{training} parameters $\param_\sigma$ ($\sigma=1,\ldots,\ell$).  
If an error estimator is available, then the parameters may be determined iteratively by the  \emph{POD-greedy} method, see~\cite{Haa17} and the references therein.

\begin{remark}
	\label{rem:blockPOD}
	MOR for coupled systems of partial or ordinary differential equations is often treated with separate basis functions for the different physical variables, as for instance, the temperature $T$ and the supply mass fraction $S$ in \eqref{eqn:wildlandfire:strong}.
	This idea of splitting the basis has been utilized in various applications, see for instance \cite{CaiIJS14,HimGB20,KerY98,LiB05,SteN13}.
	Compared to an approach that does not explicitly distinguish between different physical variables, this splitting is often observed to be advantageous in terms of preserving system properties such as stability and passivity.
\end{remark}

\subsection{Efficient approximation of nonlinear terms}
\label{subsec:onlineOfflinePOD}
To ensure the usability of the ROM~\eqref{eqn:ROM:2}, one is interested in an efficient \emph{offline/online decomposition}. 
This means that in the first stage, the \emph{offline phase}, one can spend some time constructing the ROM~\eqref{eqn:ROM:2}. 
In contrast, in the second stage, the \emph{online phase}, only a little time is available to compute an (approximate) solution for many different parameter values via the ROM~\eqref{eqn:ROM:2}. Analyzing the
ROM~\eqref{eqn:ROM:2}, we observe that the matrices $\Ared_\nu\in\R^{\ROMdim\times\ROMdim}$ ($\nu=1,\ldots,k$) can be
computed as soon as the reduced basis is determined. Hence, the matrices can be constructed in the offline phase. 
The evaluation of the linear term $\Ared(\param)\ROMstate$ thus scales with the dimension of the reduced basis and can be computed efficiently in the online phase. The situation is different for the nonlinear part, encoded in $\fred$.
Although formally, this is a mapping from $\R^{\ROMdim}\times\paramSet$ to $\R^{\ROMdim}$, its computation requires the evaluation of $\calF$, thus scaling with the dimension~$\stateDim$ of the original model, the so-called \emph{lifting bottleneck}. 
To avoid this issue, different approaches are proposed in the literature, among them \cite{RewW03,BarMNP04,AstWWB08,Gu11,ChaS10,CarFCA13,KraW19}. 

Here, we focus on the DEIM~\cite{ChaS10,DrmG16}, which is the discrete counterpart of EIM \cite{BarMNP04}. 
For a given matrix $U\in\R^{n\times m}$ of rank $m$, the DEIM approximation of $F(\state,\param)$
(cf.~\cite[Def.~3.1]{ChaS10}) on $\spann(U)$ is given by
\begin{equation}
	\label{eqn:DEIMapproximation}
	F_\DEIM(\state,\param) \vcentcolon= U(\mathbb{S}^\top U)^{-1}\mathbb{S}^\top F(\state,\param),
\end{equation}
where $\mathbb{S}\in\R^{n\times m}$ is a matrix obtained by selecting certain columns of the $n\times n$ identity matrix. 
The matrix $U$ is typically chosen by finding a best approximation of $F(\state,\param_\sigma)$ for $\sigma=1,\ldots,\ell$, i.e., by solving the POD optimization problem (respectively computing the SVD) for the data $F(\state,\param_\sigma)$. 
Note that the DEIM approximation indeed corresponds to an interpolation along the rows selected by $\mathbb{S}^\top$, since
\begin{displaymath}
	\mathbb{S}^\top F_\DEIM(\state,\param) = \mathbb{S}^\top F(\state,\param).
\end{displaymath}
Within the ROM, the approximation of the nonlinear term is thus given by
\begin{equation}
	\label{eqn:DEIMapproximation:ROM}
	\fred_\DEIM(\ROMstate,\param) \vcentcolon= V^\top F_\DEIM(V\ROMstate,\param).
\end{equation}
Let us emphasize that the computation of $\fred_\DEIM$ requires only certain rows of $F$, such that an efficient numerical implementation of $\fred_\DEIM$ that is independent of the full dimension is possible, thus avoiding the lifting bottleneck. 
For details we refer to \cite{ChaS10}. 
The question that remains to be answered is the choice of the selection matrix $\mathbb{S}$. 
While traditionally the columns of $I_n$ are selected via a greedy search~\cite{BarMNP04,ChaS10}, it was recently discovered that it is favorable to compute it by applying a QR-decomposition with column pivoting to $U^\top$, see~\cite{DrmG16} for further details. 

\subsection{Nonlinear MOR via transformation operators}
\label{subsec:shiftedPOD}

The main drawback of the MOR approach based on linear subspaces is its intrinsic separation of space and time, which cannot always be justified in applications. 
Prominent examples with strong space-time coupling are the advection equation and the wave equation. 
Given the snapshot matrices presented in \Cref{fig:snapshot}, separation of space and time implies that structures along the two principal directions can be approximated efficiently, while diagonally evolving structures are difficult to approximate by a small number of dyadic products. 
To remedy this issue, one may want to allow the modes $\varphi_i$ to evolve over time along with the solution. In \cite{BlaSU20}, this is achieved by replacing the approximation~\eqref{eqn:MORansatz} with
\begin{equation}
	\label{eqn:transformedMORansatz}
	\ROMapprox(t;\param) \vcentcolon= \sum_{i=1}^\ROMdim \ROMstate_i(t;\param)\calT_i(\Path_i(t;\param))\varphi_i,
\end{equation}
with transformation matrices $\calT_i(\Path_i)\in\R^{\stateDim\times\stateDim}$ and unknown path variables $\Path_i$. 
For the sake of presentation, we assume $\Path_i(t;\param)\in\R$.
The optimization problem~\eqref{eqn:PODminimization} with the modified ansatz~\eqref{eqn:transformedMORansatz} is then given by
\begin{equation}
	\label{eqn:shiftedPODminimization}
	\left\{\begin{aligned}
		&\min \frac{1}{2} \sum_{\sigma=1}^\ell \int_0^{\finalTime} \bigg\| \state(t;\param_\sigma) - \sum_{i=1}^\ROMdim \ROMstate_i(t;\param_{\sigma}) \calT_i(\Path_i(t;\param_{\sigma}))\varphi_i \bigg\|^2 \mathrm{d}t\\
		&\mathrm{s.t.}~\varphi_i\in\R^{\stateDim}~\text{with}~\|\varphi_i\| = 1\quad\text{for}~i=1,\ldots,\ROMdim, \text{ and}\\
		&\phantom{\mathrm{s.t.}}~\ROMstate_i(\cdot;\param_\sigma),\Path_i(\cdot;\param_\sigma)\in L^2(0,\finalTime)\quad\text{for}~\sigma=1,\ldots,\ell.
	\end{aligned}\right.
\end{equation}
A theoretical analysis of this optimization problem can be found in~\cite[sec.~4]{BlaSU20}, whereas corresponding numerical algorithms for the special case that $\calT$ is a discretized shift operator have been proposed in \cite{ReiSSM18,Rei21,SchRM19}, known as sPOD.
For further details, we refer to the forthcoming numerical examples presented in \cref{sec:numericalExample,sec:Gaussian}.
It turns out that for an efficient computation, the number of different path variables and transformation matrices should be kept as small as possible (cf.~\cite[sec.~7.1]{BlaSU20}), yielding
\begin{equation}
	\label{eqn:transformedMORansatz:2}
	\ROMapprox(t;\param) = \sum_{\rho=1}^q \sum_{i=1}^{\ROMdim_\rho} \ROMstate_{\rho,i}(t;\param)\calT_\rho(\Path_\rho(t;\param))\varphi_{\rho,i},\qquad \ROMdim \vcentcolon= \sum_{\rho=1}^q \ROMdim_\rho,
\end{equation}
i.e., multiple modes $\varphi_{\rho,i}$ are transformed by the same matrix $\calT_\rho(\Path_\rho)$. 
Although first numerical studies suggest that it is also possible to optimize over the transformation matrices $\calT_i$, see for instance \cite{RimPM20,KraSR21}, we focus here on the special case that these matrix functions are given a-priori.  
For the wildland fire application discussed above, we will use discrete approximations of the shift operator. 

\begin{example}
	To demonstrate the benefits of using transformed modes, we consider the linear wave equation with periodic boundary conditions, whose solution is presented in \Cref{fig:snapshot:waveEquation}. 
	Using POD, at least $\ROMdim=49$ modes are required to represent the solution with a relative error of less than \num{1e-5}. 
	In contrast, using the transformed modes ansatz in~\eqref{eqn:transformedMORansatz} with a shift operator as transformation, the solution can be exactly, i.e., with no approximation error, represented with $\ROMdim=2$ (respectively $\rho=2$, $\ROMdim_1 = 1$, $\ROMdim_2=1$ in~\eqref{eqn:transformedMORansatz:2}), see \cite[Ex.~4.4]{BlaSU20}.
\end{example}

If we substitute the ansatz~\eqref{eqn:transformedMORansatz} we immediately observe that $\calT_i(\Path_i)\varphi_i$ has to be differentiable with respect to $\Path_i$, which we assume in the following. 
Let us denote this derivative with $\calT_i'(\Path_i)\varphi_i$ and define the path-dependent matrices
\begin{align*}
		V(\Path) &\vcentcolon= \begin{bmatrix}
			\calT_1(\Path_1)\varphi_1 & \cdots & \calT_\ROMdim(\Path_\ROMdim)\varphi_\ROMdim
		\end{bmatrix}\in\R^{\stateDim\times\ROMdim}, & 
		W(\Path) &= \begin{bmatrix}
			\calT_1'(\Path_1)\varphi_1 & \cdots & \calT_\ROMdim'(\Path_\ROMdim)\varphi_\ROMdim
		\end{bmatrix}\in\R^{\stateDim\times \ROMdim},
\end{align*}
for the approximation ansatz~\eqref{eqn:transformedMORansatz} and
\begin{align*}
		V(\Path) &\vcentcolon= \begin{bmatrix}
			\calT_1(\Path_1)\varphi_{1,1} & \cdots & \calT_1(\Path_1)\varphi_{1,r_1} & \cdots & \calT_q(\Path_q)\varphi_{q,1} & \cdots & \calT_q(\Path_q)\varphi_{q,r_q}
		\end{bmatrix}\in\R^{\stateDim\times\ROMdim}, \\
		W(\Path) &\vcentcolon= \begin{bmatrix}
			\calT_1'(\Path_1)\varphi_{1,1} & \cdots & \calT_1'(\Path_1)\varphi_{1,r_1} & \cdots & \calT_q'(\Path_q)\varphi_{q,1} & \cdots & \calT_q'(\Path_q)\varphi_{q,r_q}
		\end{bmatrix}\in\R^{\stateDim\times \ROMdim}, 
\end{align*}
for~\eqref{eqn:transformedMORansatz:2}.  
Similarly as in the Galerkin framework, a ROM can be obtained by projecting the FOM~\eqref{eqn:FOM} with $V(\Path)$. 
To obtain equations for the path variables, we perform an additional Petrov--Galerkin projection with $V(\Path)$ and $W(\Path)D(\ROMstate)$ with $D(\ROMstate)$ defined as
\begin{equation*}
	D(\ROMstate) \vcentcolon= \diag(\ROMstate_1,\ldots,\ROMstate_{\ROMdim})\in\R^{\ROMdim\times \ROMdim}
\end{equation*}
when using the approximation ansatz \eqref{eqn:transformedMORansatz}, or defined as
\begin{equation*}
	D(\ROMstate) \vcentcolon= \blkdiag(\ROMstate_1,\ldots,\ROMstate_q)\in\R^{\ROMdim\times q}\quad\text{with}\quad \ROMstate_\rho \vcentcolon=
	\begin{bmatrix}
		\ROMstate_{\rho,1} & \cdots & \ROMstate_{\rho,\ROMdim_\rho}
	\end{bmatrix}
	^\top \text{ for }\rho=1,\ldots,q
\end{equation*}
when using the approximation ansatz \eqref{eqn:transformedMORansatz:2}.
These projection matrices naturally arise when enforcing minimization of the residual, cf.~\cite[sec.~5]{BlaSU20}.
The corresponding ROM is given as
\begin{equation}
	\label{eqn:transformedROM}
	\begin{aligned}
	\begin{bmatrix}
		I_{\ROMdim} & 0\\
		0 & D(\ROMstate)^\top
	\end{bmatrix}&\begin{bmatrix}
		M_1(\Path) & N(\Path)\\
		N(\Path)^\top & M_2(\Path)
	\end{bmatrix}\begin{bmatrix}
		I_{\ROMdim} & 0\\
		0 & D(\ROMstate)
	\end{bmatrix}\begin{bmatrix}
		\dot{\ROMstate}\\
		\dot{\Path}
	\end{bmatrix}\\
	&= \begin{bmatrix}
		\Ared_1(\Path,\param) 								& 0\\
		D(\ROMstate)^\top\Ared_2(\Path,\param) 	& 0
	\end{bmatrix}\begin{bmatrix}
		\ROMstate\\
		\Path
	\end{bmatrix} + \begin{bmatrix}
		\fred_1(\Path,\ROMstate,\param)\\
		D(\ROMstate)^\top \fred_2(\Path,\ROMstate,\param)
	\end{bmatrix}
	\end{aligned}
\end{equation}
with
\begin{align*}
	M_1(\Path) &\vcentcolon= V(\Path)^\top V(\Path)\in\R^{\ROMdim\times\ROMdim}, &
	\Ared_1(\Path,\param) &\vcentcolon= V(\Path)^\top A(\param)V(\Path)\in\R^{\ROMdim\times\ROMdim}, \\
	M_2(\Path) &\vcentcolon=  W(\Path)^\top W(\Path) \in \R^{\ROMdim\times\ROMdim},&
	\Ared_2(\Path,\param) &\vcentcolon= W(\Path)^\top A(\param)V(\Path)\in\R^{\ROMdim\times\ROMdim}, \\
	N(\Path) &\vcentcolon= V(\Path)^\top W(\Path)\in\R^{\ROMdim\times\ROMdim}, &
	\fred_1(\Path,\ROMstate,\param) &\vcentcolon= V(\Path)^\top F(V(\Path)\ROMstate,\param)\in\R^{\ROMdim},\\
	 & &
	\fred_2(\Path,\ROMstate,\param) &\vcentcolon= W(\Path)^\top F(V(\Path)\ROMstate,\param)\in\R^{\ROMdim}.
\end{align*}
Note that~\eqref{eqn:transformedROM} is a nonlinear system of equations even if the original model~\eqref{eqn:FOM} is linear. 
Let us emphasize, that in addition to a reduction in space, the ROM~\eqref{eqn:transformedROM} typically needs considerably less time steps in the temporal discretization for an accurate solution, see \cite[sec.~7.3]{BlaSU20}. 

\begin{remark}
	\label{rem:additionalPODmodes}
	If the system at hand features both transport and non-transported phenomena, it may be reasonable to add POD modes to the approximation ansatz; see for instance \cite{SchRM19}. 
	This can be achieved by adding an additional transformation matrix that equals the identity matrix and setting the corresponding path variable to zero.  
	Note that in this case, no further equation is required for the POD path variable.
\end{remark}

\section{Efficient Offline/Online Decomposition}
\label{sec:onlineOfflineDecomposition}

A quick analysis of the ROM based on transformed modes~\eqref{eqn:transformedROM} reveals that an efficient offline/online decomposition as described in \cref{subsec:onlineOfflinePOD} for standard projection-based MOR is not readily available, even if the original dynamics are linear. 
This is because the transformation operators applied to the modes have to be evaluated for different values of the path variables. Thus, in this section, we discuss further approximation steps that have to be implemented to ensure an efficient offline/online decomposition.

\begin{remark}
	\label{rem:equivariance}
	If only a single transformation matrix is used, i.e., $q=1$ in~\eqref{eqn:transformedMORansatz:2}, which is unitary and satisfies 
	\begin{equation}
		\label{eqn:relationBetweenTAndItsDerivative}
		\calT'(p)\varphi = \calT(p)P\varphi \quad \text{for all }p\in\R, \varphi\in\R^n
	\end{equation}
	for some matrix $P\in\R^{n\times n}$, and if the right-hand side of~\eqref{eqn:FOM} is equivariant with respect to this transformation, then the explicit dependency of the ROM~\eqref{eqn:transformedROM} on the transformation matrix and thus on the path variable can be avoided, see~\cite[Rem.~6.3]{BlaSU20} for a similar observation.
	In this case, an efficient offline/online decomposition can be obtained using the methods discussed in \cref{subsec:onlineOfflinePOD}. 
	An example for a unitary transformation that satisfies an infinite-dimensional analog of \eqref{eqn:relationBetweenTAndItsDerivative} is the shift operator with periodic boundary conditions, where the spatial differentiation operator assumes the role of $P$, cf.~\cite[Ex.~5.12]{BlaSU20}.
	If $q>1$, then one may still be able to show that the evaluation of~\eqref{eqn:transformedROM} can be computed without repeated evaluations of the transformation operators, see~\cite{BlaSU21a} for details. This is, however, restricted to exceptional cases and not possible in general.
\end{remark}

\subsection{Efficient approximation of path-dependent matrices}
\label{sec:innerProducts}

One significant difficulty in evaluating the ROM \eqref{eqn:transformedROM} compared to classical linear model reduction methods is that the matrices $V$ and $W$, and thus also the reduced matrices, depend on the path variable $\Path$, which is one of the ROM state variables. 
This is problematic since assembling the reduced matrices requires computations scaling with the FOM dimension, which needs to be avoided in the online phase.

To this end, we propose to precompute the reduced matrices depending on the path variables for different values of $\Path$ and construct approximants for these matrix-valued functions based on the sampled data in the offline phase. In the online phase, we can then replace $M_1$, $N$, $M_2$, $\Ared_1$, and~$\Ared_2$ by their respective approximants to avoid computations scaling with the FOM dimension. The approximants can, for instance, be constructed via interpolation or via regression. Recall that we assume the affine parameter decomposition~\eqref{eqn:affineDecomposition} for the matrix $A$, leading to the affine decompositions
\begin{equation*}
	\Ared_1(\Path,\param) = \sum_{\nu=1}^k q_\nu(\param)\Ared_{1,\nu}(\Path)\quad \text{and}\quad \Ared_2(\Path,\param) = \sum_{\nu=1}^k q_\nu(\param)\Ared_{2,\nu}(\Path),
\end{equation*}
for the reduced matrices. In particular, the strategy outlined above is applied to the matrices $\Ared_{i,\nu}$ for $i=1,2$ and $\nu=1,\ldots,k$.

Concerning the sampling procedure, we make the following observation. Each entry of the reduced matrices depends at most on two components of the path variable~$\Path(t;\param)\in\R^q$. 
If $q\leq 2$, then we can simply sample the entire matrix in $\R^q$. If, on the other hand, $q>2$, it might be advantageous to instead sample each entry of each of the matrices separately, thus alleviating the curse of dimensionality for larger $q$. 
A further reduction of the sampling space may be achieved via the method of active subspaces \cite{Con15}, see \cref{sec:MORforWildlandFire} for further details.

For some choices of the transformation $\calT$, the range of $\Path$ values is naturally bounded. This allows us to restrict the sampling to some finite domain. 
For instance, when considering the shift operator on a finite computational domain $\Omega\subseteq \R$ with periodic boundary conditions, each entry of $\Path$ can be restricted to a domain of the same length as $\Omega$.

\begin{remark}
	\label{rem:singleTransformation}
	A notable special case arises when a single transformation $\calT \vcentcolon=\calT_1 = \ldots = \calT_q$ is used, which additionally satisfies \eqref{eqn:relationBetweenTAndItsDerivative} as well as
	\begin{equation}
		\label{eqn:niceTrafos}
		\calT^\top\!(p_1)\calT(p_2)=\calT(\tau(p_1,p_2))\quad \text{for all } p_1,p_2\in \R
	\end{equation}
	for some mapping $\tau\colon \R\times\R\rightarrow \R$. 
	In this case, it is sufficient to sample the matrices on the left-hand side of~\eqref{eqn:transformedROM} in a one-dimensional space. 
	This is, for instance, the case for the shift operator on a finite domain $\Omega\subseteq \R$ with periodic boundary conditions.
	In that case, the transformation operator satisfies \eqref{eqn:niceTrafos} where $\tau$ is defined by $\tau(p_1,p_2) \vcentcolon= p_2-p_1$.
	Let us emphasize that this property is, for example, exploited in \cite[sec.~3.4]{CagMS19} for a hyper-reduction in the context of the Burgers' equation with periodic boundary conditions.
\end{remark}

\subsection{DEIM with transformation operators}
\label{subsec:shiftedDEIM}

The methodology outlined in the last subsection allows us to achieve an efficient offline/online decomposition for those ROM terms originating from linear terms in the FOM~\eqref{eqn:FOM}.
So far, we only addressed those nonlinearities that resulted from the transformation operators that render the considered approximation ansatz \eqref{eqn:transformedMORansatz} or \eqref{eqn:transformedMORansatz:2} nonlinear.
In this subsection, we focus on the second source of nonlinearity: the FOM nonlinearity $F$, which occurs in the definitions of $\fred_1$ and $\fred_2$ and prevents an efficient offline/online decomposition of those terms.
To this end, we aim to extend the idea of EIM/DEIM to our setting. For notational convenience, we restrict ourselves to the approximation ansatz \eqref{eqn:transformedMORansatz:2}, noting that~\eqref{eqn:transformedMORansatz} can be handled analogously.

The first step of DEIM, see the discussion in \cref{subsec:onlineOfflinePOD}, consists of finding a suitable linear subspace for approximating the FOM nonlinearity $F$. As outlined earlier, one popular way is to apply POD to snapshots of $F$.
However, when considering systems with transport, this transport does not only affect the snapshots of the state $\state$, but often also the snapshots of the nonlinearity (see for instance \Cref{fig:nonlinSnapshots} and \cite{SarG20}). 
As a consequence, we propose a similar approximation ansatz for the nonlinearity as for the FOM state, i{.}e{.}, we search for $\psi_{\rho,i}\in \R^n$ and $\DEIMcoeff_{\rho,i}\colon \timeInt\times\paramSet\rightarrow \R$ for $i=1,\ldots,\DEIMdim_\rho$ and $\rho=1,\ldots,q$ with $\DEIMdim\vcentcolon=\sum_{\rho=1}^q \DEIMdim_\rho\ll n$ and
\begin{subequations}
	\begin{align}
		\label{eqn:transformedNonlinearityAnsatz}
		F(\state(t;\param),\param) &\approx \sum_{\rho=1}^q \sum_{i=1}^{\DEIMdim_\rho} \DEIMcoeff_{\rho,i}(t;\param)\calT_\rho(\Path_\rho(t;\param))\psi_{\rho,i} = U(\Path(t;\param))\DEIMcoeff(t;\param)\\
		\label{eqn:transformedU}
		\text{with}\quad U(\Path) &\vcentcolon=
		\begin{bmatrix}
			\calT_1(\Path_1)\psi_{1,1} & \cdots & \calT_1(\Path_1)\psi_{1,\DEIMdim_1} & \cdots & \calT_q(\Path_q)\psi_{q,\DEIMdim_q}
		\end{bmatrix} \in\R^{n\times m},\\ 
		\DEIMcoeff(t;\mu) &\vcentcolon=
		\begin{bmatrix}
			\DEIMcoeff_{1,1}(t;\mu) & \cdots & \DEIMcoeff_{1,\DEIMdim_1}(t;\mu) & \cdots & \DEIMcoeff_{q,\DEIMdim_q}(t;\mu)
		\end{bmatrix}^\top \in \R^m.
	\end{align}
\end{subequations}
Here, we assume that the transformation matrices and paths used for approximating the state, i.e.,~$\calT_\rho(\Path_\rho)$, are also reasonable transformations for approximating the nonlinearity.
Especially for the wildland fire model, we expect this assumption to be met, since the nonlinearity is local in the sense that the nonlinearity at a certain point in the spatial domain $\Omega$ only depends on the state at the same point.
In general, there may be nonlinearities where this assumption is not satisfied and where other paths or even other transformation operators have to be used for approximating the nonlinearity.
However, such considerations are not within the scope of this paper.

The ansatz functions $\psi$ can be determined based on snapshots of $F$ in the same manner as we determine the ansatz functions $\varphi$ based on snapshots of the state $\state$ by (approximately) solving the minimization problem~\eqref{eqn:shiftedPODminimization}. 
Following the EIM/DEIM methodology, we enforce equality in~\eqref{eqn:transformedNonlinearityAnsatz} in $\DEIMdim\in\N$ selected rows.  
The selected rows within the DEIM framework are typically chosen based on the leading matrix on the right-hand side of \eqref{eqn:transformedNonlinearityAnsatz}, i.e., on the matrix $U(\Path)$. For different values of $\Path$, we may thus obtain other row selections, yielding a path-dependent selector matrix $\mathbb{S}(\Path)$.  Let us point out that in the online phase, we never explicitly need to assemble the selector matrix itself, but we only need the indices of its non-zero rows. 
The DEIM approximation with transformed modes, referred to as \emph{shifted DEIM} ($\sDEIM$), is thus given by
\begin{equation}
	\label{eqn:DEIMapproximation:transforming}
	F_\sDEIM(\Path,\state,\param) \vcentcolon= U(\Path)(\mathbb{S}(\Path)^\top U(\Path))^{-1}\mathbb{S}(\Path)^\top F(\state,\param).
\end{equation}
As before, we note that, as a consequence of the factor $\mathbb{S}(\Path)^\top$ on the right-hand side of \eqref{eqn:DEIMapproximation:transforming}, we only have to evaluate $\DEIMdim$ rows of the nonlinearity $F$.
Thus, if the Jacobian of the nonlinearity $F$ is sparse, this means that we do not need the full vector $V(\Path)\ROMstate\in\R^n$, but only a small number $\hat{\DEIMdim}(\Path)\ll n$ of its entries.
Exploiting this sparsity, we can rewrite the right-hand side of \eqref{eqn:DEIMapproximation:transforming} as
\begin{equation}
	\label{eqn:DEIMnonlinearity}
	\mathbb{S}(\Path)^\top F(\state,\param) = \widetilde{F}\left(\Path,\widetilde{\mathbb{S}}(\Path)^\top \state,\param\right).
\end{equation}
Here, $\widetilde{F}(\Path,\cdot,\param)\colon \R^{\hat{\DEIMdim}(\Path)}\rightarrow \R^{\DEIMdim}$ is the restriction of $F$ to those rows of $F$ required to evaluate the right-hand side of \eqref{eqn:DEIMapproximation:transforming} and taking only a reduced number of input arguments based on the sparsity pattern of the Jacobian of $F$.
In particular, the function $\widetilde{F}(\Path,\cdot,\param)$ can be implemented such that its evaluation does not scale with the full dimension $\stateDim$.
The details depend on the specific nonlinearity $F$. 
They are discussed in \cref{sec:MORforWildlandFire} for the nonlinearity of the wildland fire model.
Within the ROM~\eqref{eqn:transformedROM}, we thus obtain the approximations
\begin{align*}
	\fred_1(\Path,\ROMstate,\param) \approx \ROM{V}(\Path)^\top \widetilde{F}\left(\Path,\widetilde{V}(\Path)\ROMstate,\param\right)\qquad\text{and}\qquad
	\fred_2(\Path,\ROMstate,\param) \approx \ROM{W}(\Path)^\top \widetilde{F}\left(\Path,\widetilde{V}(\Path)\ROMstate,\param\right)
\end{align*}
with 
\begin{subequations}
\label{eqn:DEIMmatrices}
\begin{gather}
	\ROM{V}(\Path) \vcentcolon= V(\Path)^\top U(\Path)\left(\mathbb{S}(\Path)^\top U(\Path)\right)^{-1}\in\R^{\ROMdim\times \DEIMdim},\qquad\widetilde{V}(\Path) \vcentcolon= \widetilde{\mathbb{S}}(\Path)^\top V(\Path)\in\R^{\hat{\DEIMdim}(\Path)\times\ROMdim} ,\\
	\text{and}\qquad\ROM{W}(\Path) \vcentcolon= W(\Path)^\top U(\Path)\left(\mathbb{S}(\Path)^\top U(\Path)\right)^{-1}\in\R^{\ROMdim\times \DEIMdim}.
\end{gather}
\end{subequations}
Since these matrices depend on the path $\Path$, we propose, in the spirit of \cref{sec:innerProducts}, to first compute these matrices in the offline phase for different values of the path $\Path$, using, for instance, the QDEIM algorithm from~\cite{DrmG16}.
Afterward, we construct a corresponding interpolant depending on $\Path$, which we can then efficiently evaluate in the online phase.
At this point, we emphasize that in contrast to the construction of interpolants as outlined in \cref{sec:innerProducts}, where we never have to sample a space whose dimension is higher than two, here we have to sample a space whose dimension equals $q$, which is the number of entries in the path $\Path$. Again, the dimension may be reduced by exploiting active subspaces~\cite{Con15}. 

\begin{remark}
	If a single transformation operator with a single path is used, i.e., we are in the situation~\eqref{eqn:transformedMORansatz:2} with $q=1$, then the compressed nonlinearity snapshots may be written (assuming $\calT(0) = I_n$) as
	\begin{displaymath}
		U(\Path) = \calT(\Path)U(0).
	\end{displaymath}
	This motivates the choice $\mathbb{S}(\Path) \vcentcolon= \calT(\Path)\mathbb{S}_0$, with $\mathbb{S}_0$ chosen as the DEIM interpolation points selected for the matrix $U(0)$. We refer to \cite{RimPM20,SarG20} for similar ideas. Although in general it is possible to extend this idea to multiple transformations, i.e., $q\geq 2$, we emphasize that this may result in overlapping DEIM points, which in turn renders the matrix $\mathbb{S}(\Path)^\top U(\Path)$ singular. 
	This situation did indeed occur in our numerical experiments for the wildland fire, which is the reason we do not pursue this strategy in the following. 
\end{remark}

\section{Numerical Experiments}
\label{sec:numerics}

In this section, we present the application of the nonlinear model reduction scheme introduced in \cref{subsec:shiftedPOD,sec:onlineOfflineDecomposition} to the wildland fire model \eqref{eqn:truthModel}, simulated on a one-dimensional strip of length $\length=\si{1000\metre}$, yielding the computational domain $\Omega=[0,1000]$.
To this end, we first discuss how the general nonlinear model reduction framework is applied explicitly to the wildland fire model in \cref{sec:MORforWildlandFire}.
Afterward, we present numerical results for two different initial values in \cref{sec:numericalExample,sec:Gaussian}.
The initial value in \cref{sec:numericalExample} leads to a pure advective transport of the already formed combustion waves, whereas the case considered in \cref{sec:Gaussian} also includes their formation. In more detail, we use 
\begin{equation}
	\label{eq:GaussianIC}
	T^0(x) \vcentcolon= 1200\exp\left(-\tfrac{(x-500)^2}{200}\right)\qquad\text{and}\qquad S^0 \equiv 1
\end{equation}
as initial value in \cref{sec:Gaussian}, in agreement with the experiment in \cite[sec.~7.1]{ManBBCDKV08}. To obtain the initial condition for \cref{sec:numericalExample}, we simulate the wildland fire model with Arrhenius coefficient $\arrhenius=558.49\si{\kelvin}$, parameters listed in~\Cref{tab:coeffValues} (taken from \cite[sec.~7.1]{ManBBCDKV08}),  and initial condition~\eqref{eq:GaussianIC} for \si{700\second}. The solution at the final time is then taken as initial value for our numerical study in \cref{sec:numericalExample}.  
The time horizon is chosen so that in all numerical simulations to follow, the combustion waves never reach the boundary of the computational domain. 
For simplicity, we thus use periodic boundary conditions. 
Furthermore, as discussed in \cref{sec:wildfireModel}, we use the relative temperature $T-T_{\mathrm{a}}$ for the implementation instead of the absolute temperature $T$.

\begin{table}
	\centering
	\caption{Values for the fixed wildland fire model coefficients}
	\label{tab:coeffValues}
	\begin{tabular}{lrrrrrr}
		\toprule
		coefficient & $k$ [\si{\square\metre\per\second}] & $v$ [\si{\metre\per\second}] & $\tempRise$ [\si{\kelvin\per\second}] & $T_{\mathrm{a}}$ [\si{\kelvin}] & $\heatTransfer$ [\si{\per\kelvin}] & $\fuelRate$ [\si{\per\second}]\\
		value & $0.2136$ & $0$ & $187.93$ & $300$ & $4.8372\times 10^{-5}$ & $0.1625$ \\
		\bottomrule
	\end{tabular}
\end{table}

For the spatial discretization, we use the finite difference method with a $6$th order central finite difference stencil.
To this end, we decompose the spatial domain into $\nSpacePoints=3000$ equidistant intervals corresponding to a mesh width of~$\tfrac{1}{3}\si{\metre}$, yielding a total of $\stateDim \vcentcolon= 2\nSpacePoints = 6000$ degrees of freedom (DOF). 
The spatially discretized nonlinearity in~\eqref{eqn:truthModel} is given by
\begin{displaymath}
	f(S_h,T_h,\arrhenius) \vcentcolon= S_h\odot \hat{r}(T_h,\arrhenius),
\end{displaymath}
where $\odot$ denotes the Hadamard product and $\hat{r}$ is defined as an entry-wise application of $r$, see \eqref{eqn:reactionRate} 
The time integration is carried out using MATLAB's \texttt{ode45} function, which is an explicit time integration scheme with adaptive time stepping based on the Dormand Prince method. For the final time, we use $\finalTime = \si{1400\second}$ in \cref{sec:numericalExample} and $\finalTime = \si{2100\second}$ in \cref{sec:Gaussian}. 

In the numerical experiments, we construct parametric ROMs, which allow variation of the Arrhenius coefficient $\arrhenius$.
To this end, we use snapshots based on the values $\arrhenius\in\lbrace\si{540\kelvin,560\kelvin,580\kelvin}\rbrace$ in the offline phase and test the ROM for various $\arrhenius$ values from \SIrange{540}{580}{\kelvin}. 
For notational convenience, we refer to the POD ROM \eqref{eqn:ROM:2}, where the nonlinearity is approximated as outlined in
\cref{subsec:onlineOfflinePOD} via DEIM, as POD-DEIM ROM. Accordingly, our method resulting in the ROM~\eqref{eqn:transformedROM} with approximation of the nonlinear terms as in \cref{sec:onlineOfflineDecomposition} is abbreviated with sPOD-sDEIM ROM.

All relative errors stated in the following subsections are measured in a discretized $\Ltwo(0,\finalTime;\R^\stateDim)$ norm, where we used the trapezoidal rule for approximating the integral.

\subsection{MOR for wildfire model}\label{sec:MORforWildlandFire}

In this subsection, we discuss how we apply the general model reduction procedure outlined in \cref{subsec:shiftedPOD,sec:onlineOfflineDecomposition} to the wildland fire model \eqref{eqn:truthModel}.
Since the (relative) temperature and the supply mass fraction have different scales, we follow the principle described in \Cref{rem:blockPOD} and determine different basis functions for the temperature and for the supply mass fraction, respectively.
For each of these two quantities, we use \eqref{eqn:transformedMORansatz:2} as approximation ansatz with $q=2$ for the two traveling combustion waves as shown in \Cref{fig:snapshot:wildlandFire:temperature} and \Cref{fig:snapshot:wildlandFire:smf}.
Since the propagation of the combustion waves is reflected in both the snapshots of the temperature and the snapshots of the supply mass fraction, we take the same paths for the temperature and supply mass fraction, i.e., in total, we need to compute two paths in the online phase.
The method we use to compute the basis functions according to the approximation ansatz \eqref{eqn:transformedMORansatz:2} is outlined in \cref{subsec:separatedWavesOffline,subsec:gaussianOffline}, respectively.
For transforming the basis functions, we use the constant extrapolation shift operator, which has been applied in \cite{SchRM19} for constructing low-dimensional descriptions of reaction fronts in the context of combustion.
Let us emphasize that the supply mass fraction attains different (approximately) constant values before and after the reaction front, such that a shift operator with periodic boundary conditions would lead to unphysical supply mass fraction modes, which do not properly capture the shape of the respective reaction front.

Whenever the shift operator needs to be evaluated at points that are not multiples of the spatial mesh width, we use polynomial interpolation based on a cubic Lagrange polynomial as in \cite[sec.~7]{BlaSU20}.
This discretization of the shift operator leads in general to a non-differentiability of the mapping $\Path\rightarrow \calT(\Path)\varphi$ at multiples of the spatial grid size.
To remedy this issue, we propose to first compute the derivative of this mapping on the infinite-dimensional level and discretize afterward.
As a result, we obtain $\calT'(\Path)\varphi \vcentcolon= \calT_0(\Path)D\varphi$, where $\calT_0(\Path)$ is a discretization of the shift operator with zero extrapolation and $D$ is a sixth-order central finite difference approximation of the first derivative with one-sided finite differences at the boundary.

In the offline phase, the paths are determined by tracking the highest positive and negative value of the discrete spatial derivative within the temperature snapshots.
Here, the highest positive value corresponds to the front of the left-going combustion wave and the highest negative value to the right-going one, cf.~\Cref{fig:CombustionWavesFOMvsROM} for a depiction of the combustion waves.
Since this approach yields a step function for the paths, we partially replace the path via a linear interpolation in a post-processing step to obtain smooth functions. 
For the initial condition considered in \cref{sec:numericalExample}, the complete path trajectories are replaced by linear interpolants, which are based on the first and on the last value of the path with respect to time.
This reflects that the propagation speeds of the combustion waves are approximately constant.
The case considered in \cref{sec:Gaussian} additionally covers the formation of the wave profiles. Therefore we replace the path within the last $98.5$\% of the time interval by a linear interpolant. In contrast, the empirically determined path trajectory in the first $1.5$\% of the time interval is not a linear function of time. It is thus taken as obtained via the front-tracking approach. 

To be able to evaluate the terms $M_1$, $M_2$, $N$, $\Ared_1$, and $\Ared_2$ occurring in the ROM \eqref{eqn:transformedROM}, as well as $\ROM{V}$, $\ROM{W}$, and $\widetilde{V}$ in~\eqref{eqn:DEIMmatrices} efficiently, we use the approach presented in \cref{sec:innerProducts} as follows: In a first step, we compute an SVD of the path variables determined in the offline phase to determine an active subspace. 
We observe that it is sufficient to sample within the subspace
\begin{equation}
	\label{eqn:activeSubspace}
	\mathcal{U}_{\mathrm{a}}\vcentcolon=\{(\Path_1,\Path_2)\in \R^2 \mid \Path_1=-\Path_2\},
\end{equation}
showing that the two combustion waves propagate with the same velocity in opposite directions. 
This is due to the symmetric initial condition. Exploiting this symmetry, we only need to sample a one-dimensional domain $[0,\length]$, and based on the sampled data, we construct a piecewise constant interpolant in the offline phase, which can be evaluated efficiently in the online phase.

The nonlinearity of \eqref{eqn:truthModel} is given by 
\begin{equation}
	\label{eqn:truthNonlinearity}
	F(\state(t;\param),\param) = 
	\begin{bmatrix}
		\tempRise f(S_h,T_h,\arrhenius)\\
		-\fuelRate f(S_h,T_h,\arrhenius)
	\end{bmatrix}
	=
	\begin{bmatrix}
		\tempRise\\
		-\fuelRate
	\end{bmatrix}
	\otimes f(S_h,T_h,\arrhenius),
\end{equation}
where $\otimes$ denotes the Kronecker product.
Since the nonlinearity is given by a Kronecker product of a constant vector and a nonlinear function, we only need to approximate $f(S_h,T_h,\arrhenius)$ instead of~$F(\state,\param)$ by the $\sDEIM$ method outlined in \cref{subsec:shiftedDEIM}.
To this end, we use again the approximation ansatz \eqref{eqn:transformedMORansatz:2} for determining a low-dimensional description for $f(S_h,T_h,\arrhenius)$, which yields $U(\Path)$ in the notation of \cref{subsec:shiftedDEIM}.
The rest follows along the lines of \cref{subsec:shiftedDEIM}, with the only difference that the products $V(\Path)^\top U(\Path)$ and $W(\Path)^\top U(\Path)$ in \eqref{eqn:DEIMmatrices} have to be replaced by 
\begin{equation*}
	V(\Path)^\top \left(
	\begin{bmatrix}
		\tempRise\\
		-\fuelRate
	\end{bmatrix}
	\otimes U(\Path)\right) \quad\text{and}\quad W(\Path)^\top \left(
	\begin{bmatrix}
		\tempRise\\
		-\fuelRate
	\end{bmatrix}
	\otimes U(\Path)\right),
\end{equation*}
respectively. 

The incidence matrix corresponding to the sparsity pattern of the nonlinearity $f(S_h,T_h,\arrhenius)$ is given by
\begin{equation*}
	\begin{bmatrix}
		I_{\nSpacePoints} & I_{\nSpacePoints}
	\end{bmatrix}
	,
\end{equation*}
i.e., the $i$th entry of $f$ depends only on the $i$th and on the $(i+\nSpacePoints)$th entry of $\state$ for $i\in\lbrace 1,\ldots, \nSpacePoints\rbrace$.
As a consequence, the number of required rows of $V(\Path)$ in \eqref{eqn:DEIMnonlinearity} is $\hat{\DEIMdim}=2\DEIMdim$ and, thus, not dependent on the path.
Hence, we can achieve an efficient offline/online decomposition of the term $\widetilde{\mathbb{S}}(\Path)^\top V(\Path)$ as outlined at the end of \cref{subsec:shiftedDEIM}.
Finally, the modified nonlinearity $\widetilde{F}$ from \eqref{eqn:DEIMnonlinearity} is here simply given by
\begin{alignat*}{3}
	\widetilde{F}\colon \R^{2\DEIMdim} & \rightarrow \R^{\DEIMdim},\qquad 
	\begin{bmatrix} 
		\hat{T}_1 & \cdots & \hat{T}_{\DEIMdim} & \hat{S}_1 & \cdots & \hat{S}_{\DEIMdim}
	\end{bmatrix}
	^\top 
	& \mapsto 
	\begin{bmatrix}
		\hat{S}_1\hat{r}(\hat{T}_1,\arrhenius) & \cdots & \hat{S}_\DEIMdim\hat{r}(\hat{T}_\DEIMdim,\arrhenius)
	\end{bmatrix}
	^\top.
\end{alignat*}
Especially, we note that for the specific nonlinearity of the wildland fire model, $\widetilde{F}$ does not depend on the path.
Moreover, the mapping $\widetilde{F}$ allows to vary the parameter $\arrhenius$ in the ROM without requiring a parameter-affine structure in $f$.

\subsection{Initial condition with traveling wave solution}
\label{sec:numericalExample}

In this section, we discuss the numerical solution of the ROM \eqref{eqn:transformedROM} for the initial condition where the initial ignition process has
already taken place, and combustion waves have already developed. 

\subsubsection{Determination of modes}
\label{subsec:separatedWavesOffline}
\begin{figure}
	\centering
	\begin{subfigure}[t]{.33\linewidth}
		\begin{tikzpicture}
			\pgfplotsset{width=6cm,height=6cm}
			\begin{groupplot}[group style={group size=1 by 1, horizontal sep=2.5cm}]
				\nextgroupplot[enlargelimits=false,axis on top=false,xlabel=$x\vphantom{-\Path_1(t)}$,ylabel=$t$,xtick=\empty,ytick=\empty,ylabel style={rotate=-90}]
				\addplot graphics [xmin=0,xmax=1,ymin=0,ymax=1]{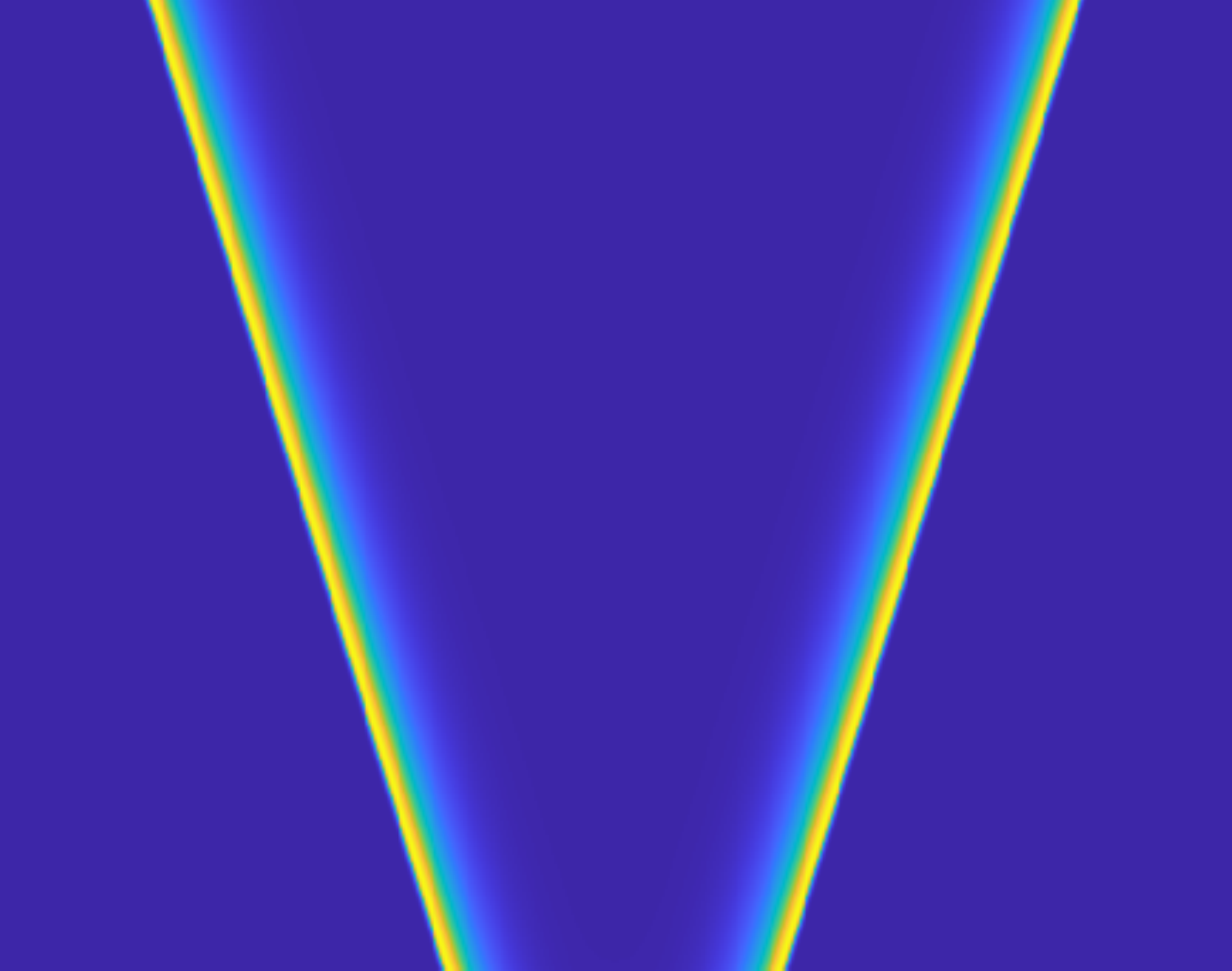};
			\end{groupplot}
			\draw[white, dashed, thick] (2.2,0) to (2.2,4.5);
		\end{tikzpicture}
		\caption{Original snapshot data}
		\label{fig:snapshot:separatedWavesAll}
	\end{subfigure}\hfill
	\begin{subfigure}[t]{.33\linewidth}
		\begin{tikzpicture}
			\pgfplotsset{width=6cm,height=6cm}
			\begin{groupplot}[group style={group size=1 by 1, horizontal sep=2.5cm}]
				\nextgroupplot[enlargelimits=false,axis on top=false,xlabel=$x\vphantom{-\Path_1(t)}$,xtick=\empty,ytick=\empty,ylabel style={rotate=-90},ylabel=$t$]
				\addplot graphics [xmin=0,xmax=1,ymin=0,ymax=1]{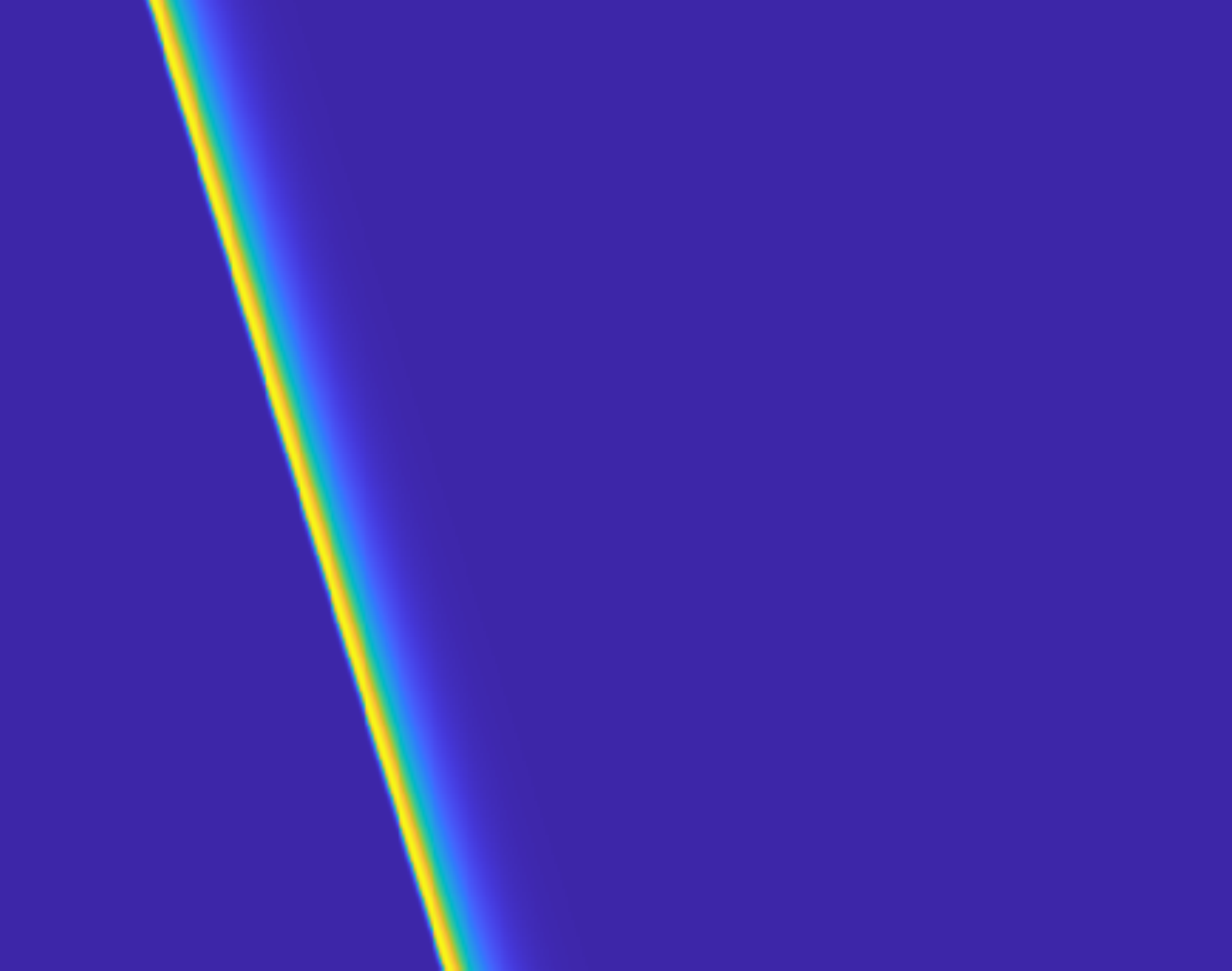};
			\end{groupplot}
		\end{tikzpicture}
		\caption{Left-going combustion wave}	
		\label{fig:snapshot:separatedWavesLeftWave}
	\end{subfigure}\hfill
	\begin{subfigure}[t]{.33\linewidth}
		\begin{tikzpicture}
			\pgfplotsset{width=6cm,height=6cm}
			\begin{groupplot}[group style={group size=1 by 1, horizontal sep=2.5cm}]
				\nextgroupplot[enlargelimits=false,axis on top=false,xlabel=$x-\Path_1(t)$,xtick=\empty,ytick=\empty,ylabel style={rotate=-90},ylabel=$t$]
				\addplot graphics [xmin=0,xmax=1,ymin=0,ymax=1]{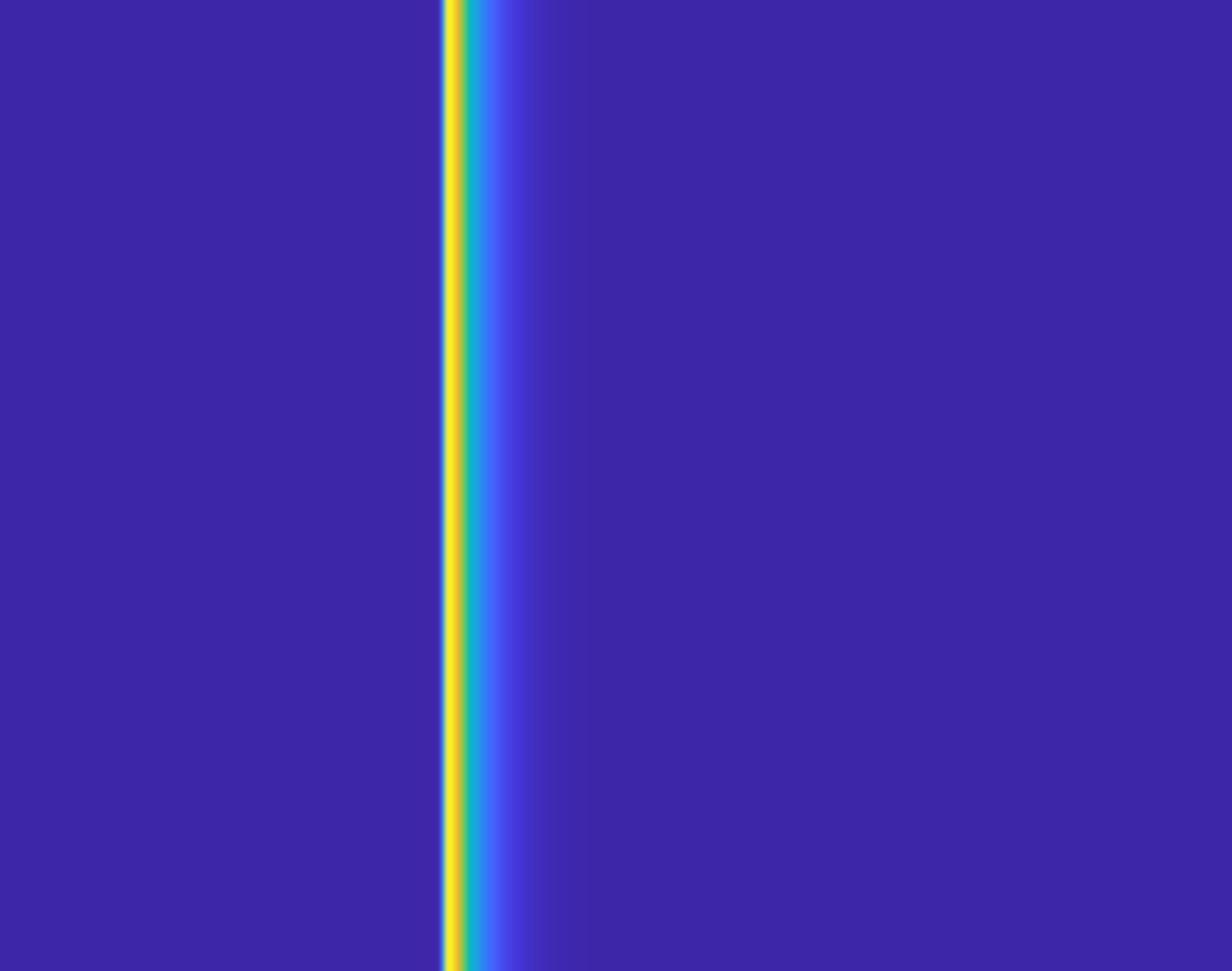};
			\end{groupplot}
		\end{tikzpicture}
		\caption{Left-going combustion wave in the co-moving frame}	
		\label{fig:snapshot:separatedWavesRightWave}
	\end{subfigure}
	\caption{Illustration of the separation procedure for the snapshot data}
	\label{fig:snapshot:separatedWaves}
\end{figure}

We use the following procedure to approximately solve the optimization problem~\eqref{eqn:shiftedPODminimization}: We first separate the two combustion waves by dividing the snapshot matrix in the middle (cf.~\Cref{fig:snapshot:separatedWavesAll}) and replacing the left/right combustion wave with zeros, see~\Cref{fig:snapshot:separatedWavesLeftWave}.  
Afterward, each set of snapshot data is now transformed to the co-moving frame, which is obtained by shifting the spatial domain according to the pre-determined path. The
resulting transformed snapshot data is low-rank, and suitable modes are now determined from the transformed snapshot
data via POD, respectively SVD. 
The approximations for each co-moving frame are then shifted back and added back together, to get an approximation of the combined snapshot data. 
The corresponding offline errors and online errors for Arrhenius coefficient $\arrhenius = 558.49\si{\kelvin}$ for different numbers of modes are presented in \Cref{tab:separatedWavesOfflineOfflineErr}, showing that this strategy seems indeed reasonable to tackle the optimization problem~\eqref{eqn:shiftedPODminimization}. 
Note that the DOF are given by the total number of transformed modes plus the additional two path variables, which we indicate with the notation in \Cref{tab:separatedWavesOfflineOfflineErr}.  
The first two modes for the left-going combustion wave are presented in \Cref{fig:separatedWavesModesLeftWave}. 
We see, in agreement with the offline errors presented in \Cref{tab:separatedWavesOfflineOfflineErr}, that already the first mode captures the dominant dynamics. 
Let us emphasize, that the coefficient for the first temperature mode varies less than \num{2e-3} percent, while the coefficient for the first supply mass fraction mode varies less than \num{9e-3} percent.
\begin{figure}
	\begin{subfigure}[t]{.49\linewidth}
		\input{ROM-separatedWaves-modes-var1}
		\caption{Relative temperature modes}
	\end{subfigure}\hfill
	\begin{subfigure}[t]{.49\linewidth}
		\input{ROM-separatedWaves-modes-var2}
		\caption{Supply mass fraction modes}
	\end{subfigure}
	\caption{The first two (normalized) modes for the left combustion wave for the relative temperature and the supply mass fraction}
	\label{fig:separatedWavesModesLeftWave}
\end{figure}
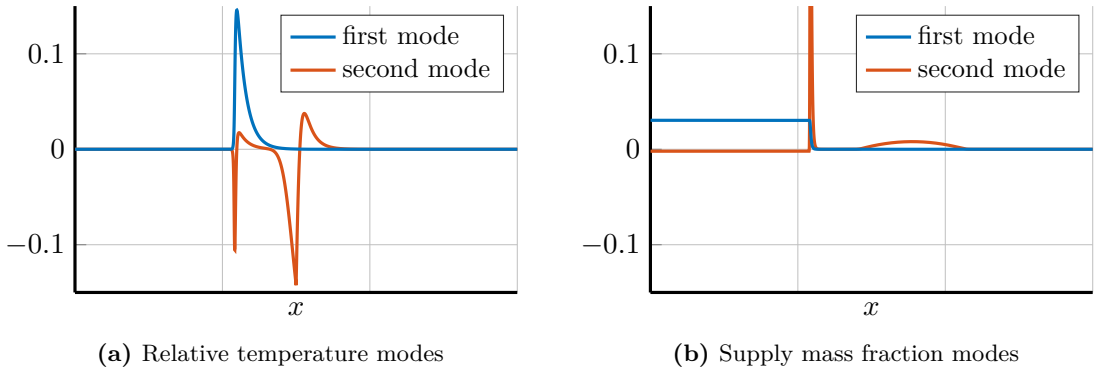
\begin{table}
	\centering
	\caption{Comparison of the relative offline and online errors for the sPOD-sDEIM ROM, trained and tested with $\arrhenius=558.49\si{\kelvin}$, for different numbers of modes}
	\label{tab:separatedWavesOfflineOfflineErr}
	\begin{tabular}{lcc@{\hspace{2em}}cc@{\hspace{2em}}r}
		\toprule
		\multirow{2}[2]{*}{DOF} & \multicolumn{2}{c}{relative offline error} & \multicolumn{2}{c}{relative online error} & \multirow{2}[2]{*}{speedup}\\\cmidrule(lr{2em}){2-3}\cmidrule(r{2em}){4-5}
		 & $T-T\mathrm{a}$ & $S$ & $T-T\mathrm{a}$ & $S$\\ \midrule
		$4+2$ & \num{1.901984e-03} & \num{4.278072e-04} & \num{6.177339e-03} & \num{1.410139e-03} & \num{130.717824}\\
		$8+2$ & \num{2.903790e-04} & \num{1.054369e-05} & \num{4.428441e-04} & \num{9.134173e-05} & \num{108.964590}\\
		$12+2$ & \num{2.274642e-04} & \num{9.216271e-06} & \num{2.433053e-04} & \num{4.919068e-05} & \num{69.599759}\\
		$16+2$ & \num{1.503261e-04} & \num{8.411403e-06} & \num{2.584320e-04} & \num{5.126845e-05} & \num{58.169859}\\
		$20+2$ & \num{1.187341e-04} & \num{7.931601e-06} & \num{1.904002e-04} & \num{4.395058e-05} & \num{43.801491}\\
		$24+2$ & \num{9.894865e-05} & \num{7.645134e-06} & \num{2.829387e-04} & \num{6.260659e-05} & \num{32.121681}\\\bottomrule
	\end{tabular}
\end{table}

\subsubsection{Efficient offline/online decomposition}
\label{subsec:separatedWaves:offlineOnline}
\begin{figure}
	\centering
	\begin{subfigure}[t]{.45\linewidth}
		\centering
		\begin{tikzpicture}
			\pgfplotsset{width=7cm,height=7cm}
			\begin{groupplot}[group style={group size=1 by 1, horizontal sep=2.5cm}]
				\nextgroupplot[enlargelimits=false,axis on top=false,xlabel=$x$,ylabel=$t$,xtick=\empty,ytick=\empty,ylabel style={rotate=-90}]
				\addplot graphics [xmin=0,xmax=1,ymin=0,ymax=1]{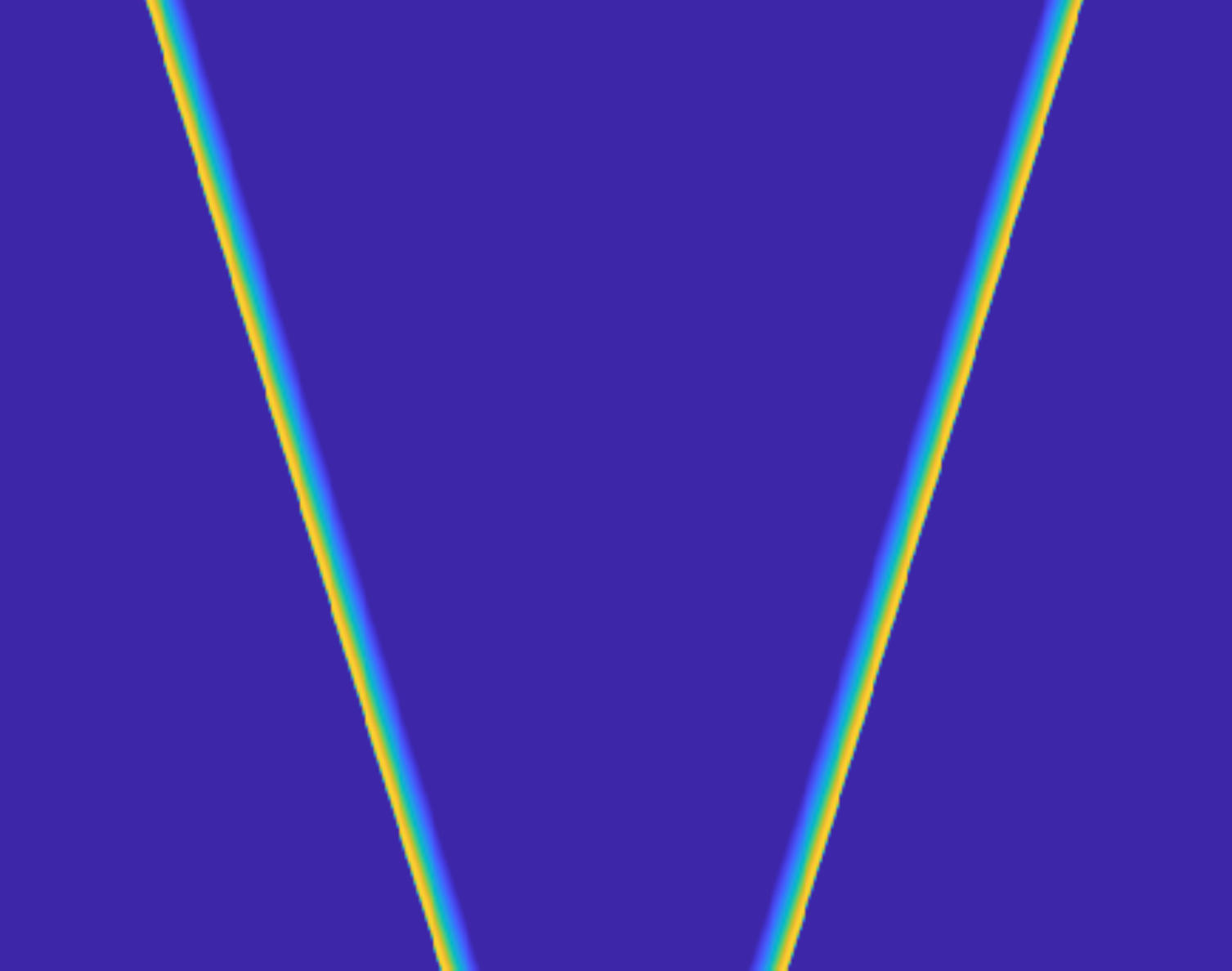};
			\end{groupplot}
		\end{tikzpicture}
		\phantom{t} 
		\caption{Separated waves initial condition}
		\label{fig:nonlinSnapshots:separatedWaves}
	\end{subfigure}\hfill
	\begin{subfigure}[t]{.53\linewidth}
		\centering
		\begin{tikzpicture}
			\pgfplotsset{width=7cm,height=7cm}
			\begin{groupplot}[group style={group size=1 by 1, horizontal sep=2.5cm}]
				\nextgroupplot[enlargelimits=false,axis on top=false,xlabel=$x$,xtick=\empty,ytick=\empty,ylabel style={rotate=-90},ylabel=$t$,colorbar,point meta min = 6.279e-6,point meta max =0.6279,colorbar style={at={(1.1,1)},anchor=north west,ymode=log}]
				\addplot graphics [xmin=0,xmax=1,ymin=0,ymax=1]{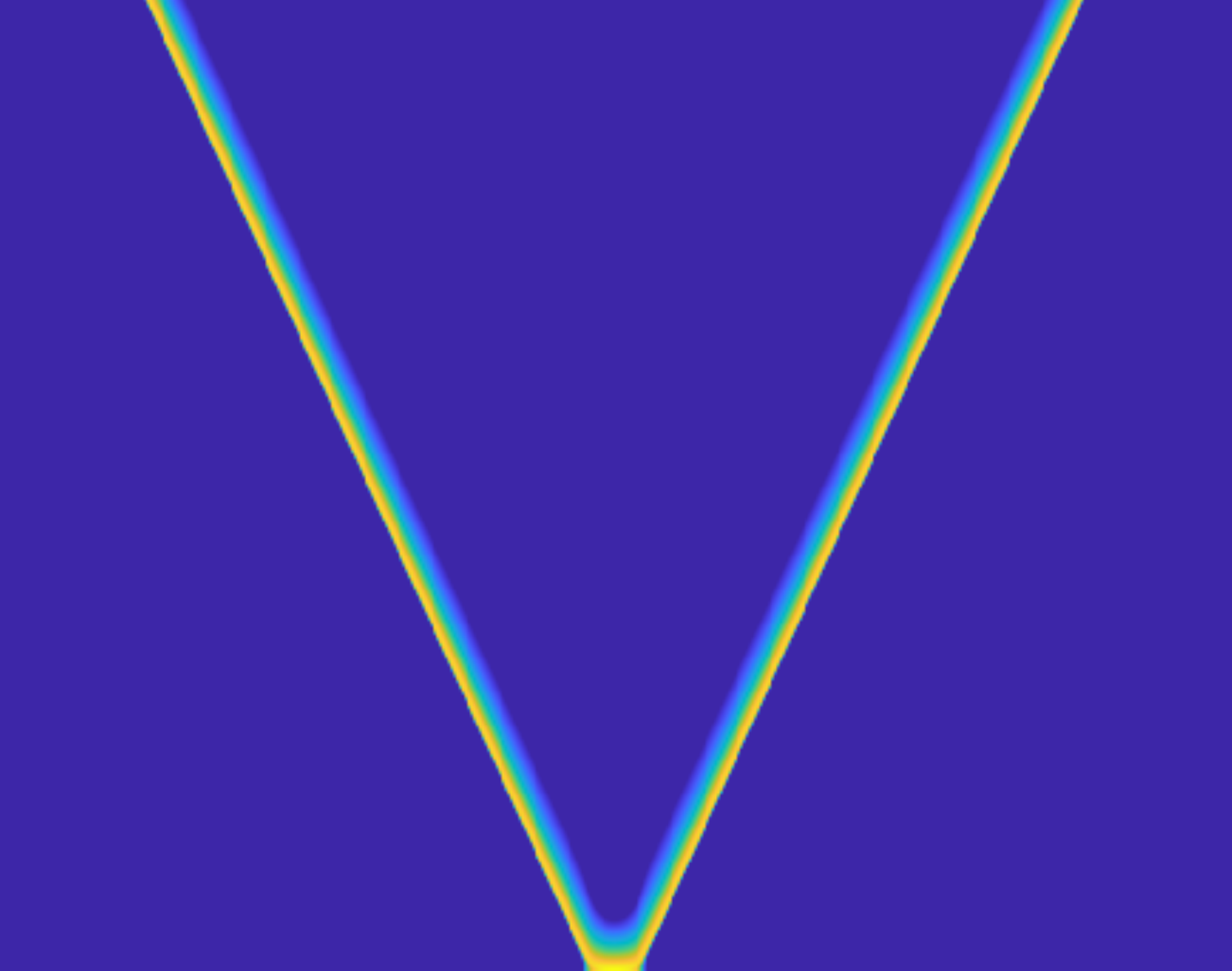};
			\end{groupplot}
		\end{tikzpicture}
		\caption{Gaussian initial condition}	
		\label{fig:nonlinSnapshots:Gaussian}
	\end{subfigure}
	\caption{Snapshots of the nonlinearity}
	\label{fig:nonlinSnapshots}
\end{figure}

For the $\sDEIM$ method as described in \cref{subsec:shiftedDEIM}, we first notice that the snapshots of the nonlinearity indeed feature similar combustion waves with the same wave speeds, see \Cref{fig:nonlinSnapshots} for an illustration.  
This justifies that we take the same transformations and the same paths used for the approximation of the state to approximate the nonlinearity.
Moreover, the mode numbers are chosen such that the number of nonlinearity modes is twice the number of modes used for the temperature and supply mass fraction (both temperature and supply mass fraction have the same number of modes).

As outlined in \cref{sec:innerProducts,subsec:shiftedDEIM}, we need to sample a couple of path-dependent matrix functions to achieve an efficient offline/online decomposition.
Therefore, as mentioned before, we make use of the active subspace defined in \eqref{eqn:activeSubspace}, and only sample within a finite subset of this subspace.
For the settings considered in this subsection, it appears to be sufficient to use $\{(\Path_1,\Path_2)\in \mathcal{U}_{\mathrm{a}} \mid \Path_1\in[0,300]\}$ as sampling domain, since the combustion waves never travel more than $300$ meters for the considered range of Arrhenius coefficients $\beta$ and for the considered time horizon.
For the discretization of the sampling domain, we use an equidistant sampling grid with a grid size of $\tfrac{20}{3}$ meters, since this leads to reasonable approximation errors and computation times for all considered mode numbers.

\subsubsection{ROM simulations}
Using the efficient offline/online decomposition, we can now compare the efficiency of the ROMs for different numbers of modes. 
To this end, we always use the same number of modes per path or reference frame, and we use the same number of modes for the temperature and for the supply mass fraction.
The total number of DOF for the ROM obtained by our approach is then given as $\ROMdim+2$, where~$\ROMdim$ is the total number of transformed modes and~$2$ is the number of frames for the considered test case, i.e., the number of path variables.
The details of the errors and speedups for different mode numbers are presented in \Cref{tab:separatedWavesOfflineOfflineErr}, showing that the relative online errors are within one order of magnitude to the offline errors. We notice that while the error in the offline phase decays with an increasing number of modes, this is not necessarily true in the online phase. Let us emphasize that the online error is subject to several different approximations, such that we cannot provide a decisive explanation for the observed behavior. A possible explanation for stagnation in the online error seems to be that with an increasing number of modes, the ROM is more sensitive to errors stemming from the spatial discretization. 
A detailed error analysis is not within the scope of this paper and subject to future research.

\begin{figure}
	\centering
	\input{separatedWaves-POD-vs-SPOD.tex}
	\caption{Comparison of computation time and accuracy for $\arrhenius = 558.49\si{\kelvin}$ for different mode numbers: sPOD-sDEIM ROM (blue stars) vs.~POD-DEIM ROM (red circles)}
	\label{fig:comp_time_vs_accuracy}
\end{figure}
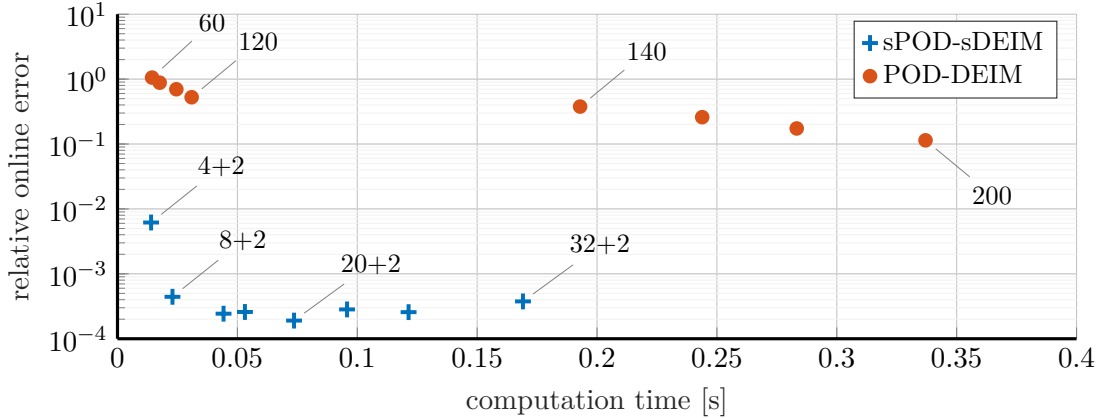

We compare POD-DEIM ROMs and sPOD-sDEIM ROMs for
the Arrhenius coefficient $\arrhenius=558.49\si{\kelvin}$ with respect to computation time and accuracy, for different numbers of
modes. To avoid irregularities in the computation time, we average the computation time over three simulations per
distinct mode number. The results are presented in \Cref{fig:comp_time_vs_accuracy}. We find that the sPOD-sDEIM ROM 
outperforms the POD-DEIM ROM, both with respect to accuracy and computation time. Note that the POD-DEIM ROM with 200 modes still results in an error that is one order of magnitude worse than all sPOD-sDEIM ROMs. 

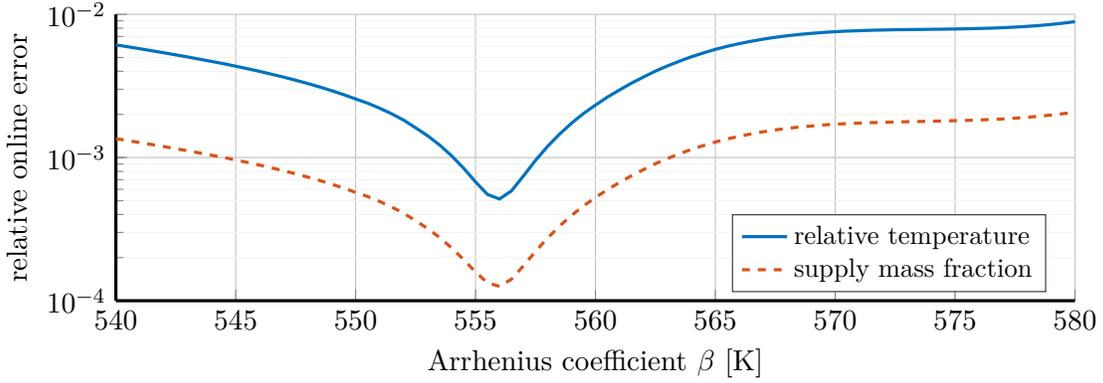
\begin{figure}
	\centering
	\input{separatedWavesParameterSamplingN3000r3-3hyper}
	\caption{Relative online error of the sPOD-sDEIM ROM with $12+2$ DOF for $\arrhenius\in\{540\si{\kelvin}, 560\si{\kelvin}, 580\si{\kelvin}\}$ over the Arrhenius coefficient $\arrhenius$ ranging from \si{540\kelvin} to \si{580\kelvin}}
	\label{fig:parameterSampling_sepWaves}
\end{figure}

To investigate the ability of the sPOD-sDEIM ROM to make predictions when evaluated at parameter values not in the training set, we use the sPOD-sDEIM ROM with $3$ modes per variable and frame, yielding a total of $3\times 2\times 2 +2=14$ DOF (including the path variable for each frame).
As training parameter values, we use three samples of the FOM with $\arrhenius\in\{\si{540\kelvin}, \si{560\kelvin}, \si{580\kelvin}\}$. We test the sPOD-sDEIM ROM for various $\arrhenius$ values within the range from \si{540\kelvin} to \si{580\kelvin}.
The resulting relative online errors for the temperature and for the supply mass fraction are depicted in \Cref{fig:parameterSampling_sepWaves}.
For computing the errors, we also simulated the FOM for these parameter values. We observe that with the settings at hand, the ROM is able to make accurate predictions over the complete sampling interval, with worst-case error less than $\si{1\percent}$. In average, our computations show that the ROM is more than $75$ times faster than the FOM with a mean error of less than 0.5\si{\percent}.

\subsection{Initial condition with combined traveling and non-transported effects}
\label{sec:Gaussian}

\begin{figure}
	\centering
	\begin{subfigure}[t]{.42\linewidth}
		\centering
		\begin{tikzpicture}
			\pgfplotsset{width=7cm,height=7cm}
			\begin{groupplot}[group style={group size=1 by 1, horizontal sep=2.5cm}]
				\nextgroupplot[enlargelimits=false,axis on top=false,xlabel=$x$,ylabel=$t$,xtick=\empty,ytick=\empty,ylabel style={rotate=-90}]
				\addplot graphics [xmin=0,xmax=1,ymin=0,ymax=1]{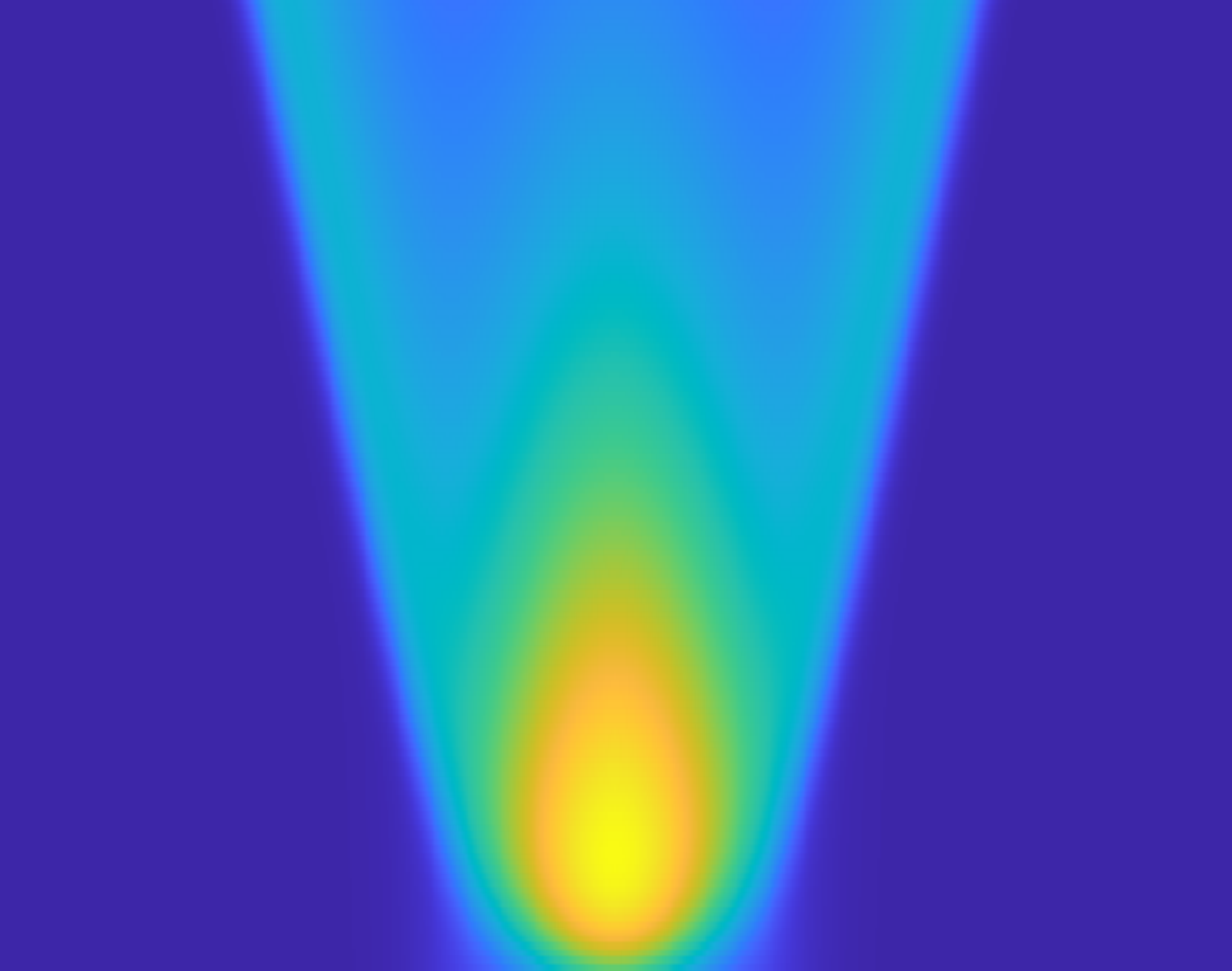};
			\end{groupplot}
			\draw[white, dotted, thick] (0,3.6) to (5.4,3.6);
			\node[white] at (0.55,1.46) {\small\beforeSwitchTag};
			\node[white] at (0.55,4) {\small\beforeSPODTag};
		\end{tikzpicture}
		\phantom{t} 
		\caption{Zoom time interval $[0,150]$}
		\label{fig:GaussianTemperatureZoomIC}
	\end{subfigure}\hfill
	\begin{subfigure}[t]{.56\linewidth}
		\centering
		\begin{tikzpicture}
			\pgfplotsset{width=7cm,height=7cm}
			\begin{groupplot}[group style={group size=1 by 1, horizontal sep=2.5cm}]
				\nextgroupplot[enlargelimits=false,axis on top=false,xlabel=$x$,xtick=\empty,ytick=\empty,ylabel style={rotate=-90},ylabel=$t$,colorbar,point meta min = 7.3e-10,point meta max =1.9e3,colorbar style={at={(1.1,1)},anchor=north west}]
				\addplot graphics [xmin=0,xmax=1,ymin=0,ymax=1]{FOM_T_Gaussian};
			\end{groupplot}
			\draw[white, dotted, thick] (0,0.26) to (5.4,0.26);
			\draw[white, dashed, thick] (0,2.1) to (5.4,2.1);
			\node[white] at (0.55,0.1) {\small\beforeSwitchTag};
			\node[white] at (0.55,0.7) {\small\beforeSPODTag};
			\node[white] at (0.55,2.83) {\small\sPODTag};
		\end{tikzpicture}
		\caption{Snapshot data}	
		\label{fig:snapshot:wildlandFire:temperature:2}
	\end{subfigure}
	\caption{Evolution of the relative temperature for the wildland fire model with Gaussian initial condition~\eqref{eq:GaussianIC}}
\end{figure}
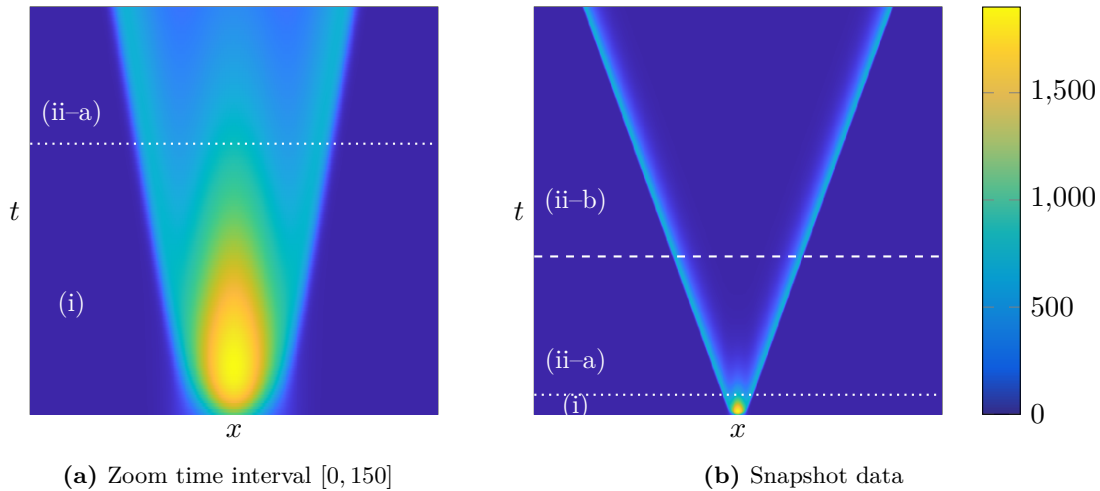

We now turn to the more challenging scenario with the initial value given in \eqref{eq:GaussianIC}. The corresponding solution of the discretized wildland fire model~\eqref{eqn:wildlandfire:strong} is shown in \Cref{fig:snapshot:wildlandFire:temperature,fig:snapshot:wildlandFire:smf}. Let us emphasize that in the beginning of the evolution of the temperature variable (cf.~\Cref{fig:GaussianTemperatureZoomIC}), a separation of the traveling waves and the ignition is not straightforward. Instead, we follow the strategy from \cite{DihDH11} and partition the time interval to construct two different surrogate models:
\begin{itemize}
	\item In the beginning of the time interval, i.e., in the area \beforeSwitchTag in \Cref{fig:GaussianTemperatureZoomIC}, we approximate the wildland fire model with standard POD-Galerkin, as described in \cref{subsec:onlineOfflinePOD}.
	\item As soon as a clear separation of the waves and the ignition is possible, i.e., in areas \beforeSPODTag and \sPODTag in \Cref{fig:snapshot:wildlandFire:temperature:2}, we use our ROM with transformation operators to capture the two wave fronts and add additional POD modes, as described in \Cref{rem:additionalPODmodes}, for the remaining dynamics.
\end{itemize}

\subsubsection{Determination of modes}
\label{subsec:gaussianOffline}
As outlined earlier, we use a POD approximation in the time interval $t\in[0,\switchTime]$, which corresponds to the area \beforeSwitchTag in \Cref{fig:snapshot:wildlandFire:temperature:2}, and a combined approximation of transformed modes and POD modes in the interval $t\in[\switchTime,\finalTime]$ with $\switchTime = \si{100\second}$, corresponding to areas \beforeSPODTag and \sPODTag. 

\begin{figure}
	\centering
	\input{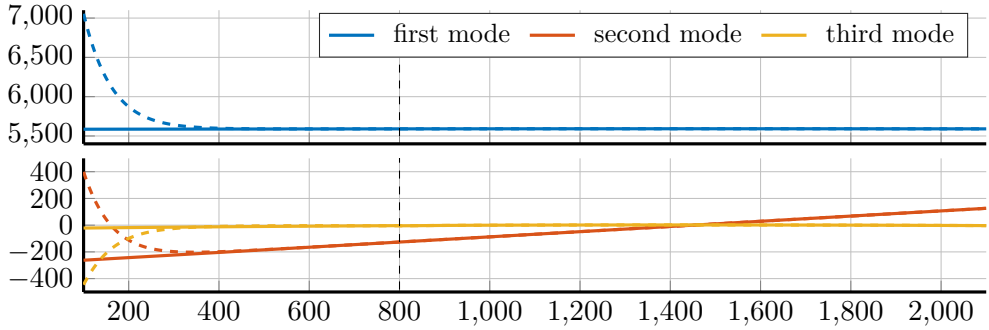}
	\caption{Extrapolated (solid) and projected (dashed) coefficients of the shifted temperature modes for the left-traveling wave}
	\label{fig:extrapolatedVsProjectedCoefficients}
\end{figure}

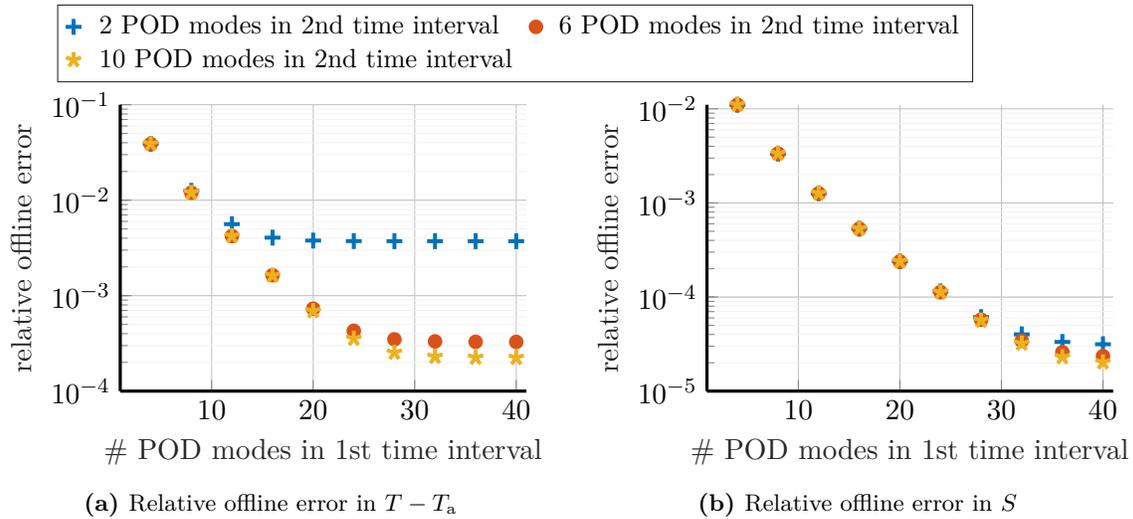
\begin{figure}
	\centering
	\begin{subfigure}[t]{.48\linewidth}
		\input{varyNumbersPODModes_T_offline}
		\caption{Relative offline error in $T-T_{\mathrm{a}}$}
		\label{fig:varyNumbersPODModes_T}
	\end{subfigure}\hfill
	\begin{subfigure}[t]{.48\linewidth}
		\input{varyNumbersPODModes_S_offline}
		\caption{Relative offline error in $S$}
		\label{fig:varyNumbersPODModes_S}
	\end{subfigure}
	\caption{Relative offline errors for different POD mode numbers in the first and in the second time interval}
	\label{fig:varyNumbersPODModes}
\end{figure}

We have to separate the traveling waves and the non-transported effects for the combined approach with transformed modes and POD modes. To minimize the impact of the non-transported effects on the modes for the traveling waves, we further separate the area after the switch and extract the transformed modes solely from area \sPODTag, which consists of $\si{65\percent}$ of the time domain after the switch for the relative temperature and $\si{80\percent}$ for the supply mass fraction. We, therefore, propose the following strategy:
\begin{enumerate}
	\item In area \sPODTag in \Cref{fig:snapshot:wildlandFire:temperature:2}, we proceed with the heuristic discussed in \cref{subsec:separatedWavesOffline}, i.e., we separate the two combustion waves by dividing the computational domain in the middle, replacing the missing parts with zeros, and then shift accordingly. For an illustration, we refer to \Cref{fig:snapshot:separatedWaves}. In the shifted frame, i.e., in \Cref{fig:snapshot:separatedWavesLeftWave}, we apply standard POD to obtain the transformed modes. 
	\item To eliminate the traveling waves also from the area \beforeSPODTag, we fit a low-degree polynomial to the coefficients $\ROMstate_{\rho,i}(t;\param)$ in \eqref{eqn:transformedMORansatz:2} for the traveling wave modes obtained within the shifted frames, and use these polynomials to extrapolate the coefficients. 
	\item We subtract the traveling wave approximations from the snapshot data and apply standard POD to capture the non-transported effects.
\end{enumerate}
\noindent
We conclude the offline phase with several important remarks:
\begin{itemize}
	\item The above strategy (and similarly the method in \cref{subsec:separatedWavesOffline}) will, in general, not provide a minimizer for the optimization problem~\eqref{eqn:shiftedPODminimization}.
		Nevertheless, these heuristic methods can be computed efficiently and deliver satisfactory results.
		In contrast, techniques that are based on solving the minimization problem~\eqref{eqn:shiftedPODminimization} or related minimization problems, cf.~\cite{Rei21,SchRM19}, involve expensive iterative solvers for an optimization problem, whose number of parameters scales with the dimension of the spatial and temporal discretization.
	\item Note that we extrapolate the coefficients instead of computing them via projection. 
		The main reason for this is that if we simply project, then the non-transported effects are, in parts, also approximated by the transformed modes. 
		This can be seen from the coefficients of the first three modes for the left-traveling temperature combustion wave depicted in \Cref{fig:extrapolatedVsProjectedCoefficients}, where the extrapolated coefficients are depicted with solid lines and the coefficients obtained via projection with dashed lines. 
		Furthermore, in the numerical experiments, we observed that choosing the extrapolation-based approach instead of projection leads to significantly smaller offline and online errors.
	\item In \Cref{fig:varyNumbersPODModes}, the offline errors are depicted for different choices for the number of POD modes in the first and in the second time interval.
		The results are based on 3 transformed modes per variable and frame and Arrhenius coefficient $\arrhenius = 558.49\si{\kelvin}$.
		We observe different error decays for the temperature and the supply mass fraction: The error decay in the temperature suggests that the numbers of POD modes chosen for the first and second time intervals should be adequately balanced.
		For instance, if only one POD mode per variable is chosen for the second time interval, it does not pay off to increase the number of used POD modes for the first time interval to values higher than $10$, since for higher values, the errors stagnate. 
		In contrast, the error decay for the supply mass fraction seems to indicate that the error is more or less independent from the number of POD modes used in the second time interval.
		These different behaviors in the temperature and supply mass fraction error are likely because at the beginning of the second time interval, there is still a fading temperature peak in the middle of the computational domain.
		However, this peak is not reflected in the supply mass fraction, which is already almost entirely consumed in the middle of the computational domain at the beginning of the second time interval.
		Consequently, adding POD modes in the second time interval is especially valuable for approximating the temperature but less significant for the supply mass fraction.
		Despite this observation, we decide, for simplicity, to always take the same number of modes for the temperature and for the supply mass fraction while noting that there is some unexploited potential for further improvement. 
\end{itemize}

\subsubsection{Efficient offline/online decomposition}
\label{subsub:gaussianOfflineOnline}

The separation of the time interval into two subintervals and the corresponding usage of two different ROMs also must be taken into account when ensuring an efficient offline/online decomposition.
For the first time interval, we employ a standard POD-DEIM ROM. We use the same number of modes for the temperature and for the supply mass fraction, whereas we use twice as many modes for approximating the nonlinearity $f$.

For the second time interval, we use a ROM with transformed modes to capture the traveling combustion waves, enriched with additional POD modes for describing the fading combustion in the middle of the computational domain.
Again, we use the same number of modes for the temperature and for the supply mass fraction, whereas we choose twice as many POD and transformed modes for the nonlinearity.
For sampling the path-dependent matrices, we exploit the active subspace as in \cref{subsec:separatedWaves:offlineOnline} and sample within the range $\{(\Path_1,\Path_2)\in \mathcal{U}_{\mathrm{a}} \mid \Path_1\in[0,500]\}$ .
This sampling domain is discretized using an equidistant grid with a grid size of $\frac{1}{3}$ meters.

In between the two time intervals, we need to switch between the two ROMs. This requires computing a new initial value for the ROM used in the second time interval.
To this end, we first compute the approximation of the full-dimensional state based on the reduced state calculated at the last time point of the first time interval.
Afterward, we use this approximation to determine the initial paths of the second time interval by using the same tracking procedure as in the offline phase, see \cref{sec:MORforWildlandFire}.
Finally, the remaining coefficients of the initial condition are computed by an orthogonal projection of the approximation onto the span of the transformed modes and the POD modes used for the second time interval.

\subsubsection{ROM simulations}
\label{subsec:Gaussian:ROMSimulations}

\begin{figure}
	\centering
	\input{Gaussian_POD_vs_SPOD.tex}
	\caption{Comparison of computation time and accuracy for  $\arrhenius = 558.49\si{\kelvin}$ with the initial condition~\eqref{eq:GaussianIC} for different mode numbers: sPOD-sDEIM ROM (blue stars) vs.~POD-DEIM ROM (red circles). 
	For the POD-DEIM ROMs, the number indicates the total number of POD modes used for this ROM. 
	For the sPOD-sDEIM ROMs the first number denotes the POD modes before the switch. 
	The numbers within the parentheses are the number of transformed modes, the number of POD modes after the switch, and the path variables. }
	\label{fig:comp_time_vs_accuracy_gaussian}
\end{figure}
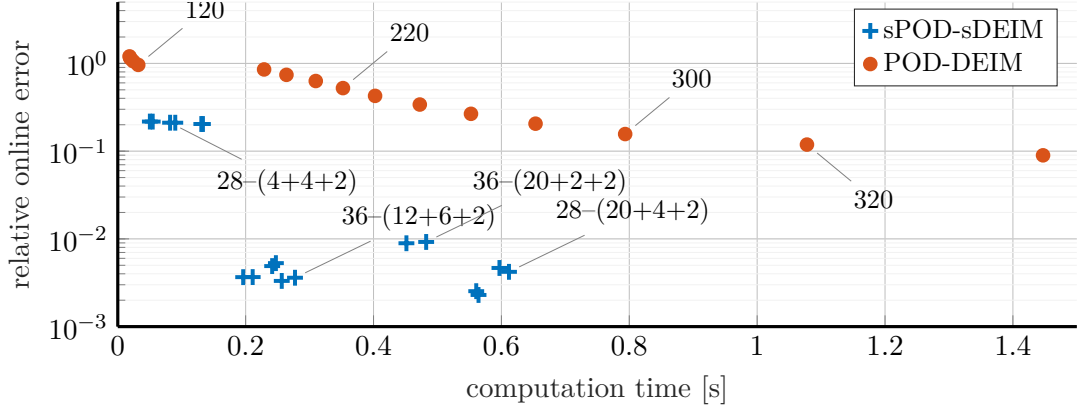

We conclude our numerical study by detailing the online performance indicators of the sPOD-sDEIM ROM.  We start by comparing a standard POD-DEIM ROM with our method in terms of computation time and accuracy. The results are presented in \Cref{fig:comp_time_vs_accuracy_gaussian}. As before, the computation times are averages over three simulations.
As in the previous section, our method outperforms standard POD with a better accuracy to computation time ratio for all considered DOF. 

\begin{figure}
	\centering
	\input{ROM-Gaussian-parameterTesting}
	\caption{Relative online error of the sPOD-sDEIM ROM with $28$ POD modes before switching and $12$ transformed modes,  $4$ POD modes, and $2$ path variables after switching, trained for $\arrhenius\in\{540\si{\kelvin}, 560\si{\kelvin}, 580\si{\kelvin}\}$, over the Arrhenius coefficient $\arrhenius$ ranging from \si{540\kelvin} to \si{580\kelvin}}
	\label{fig:parameterSampling_Gaussian}
\end{figure}
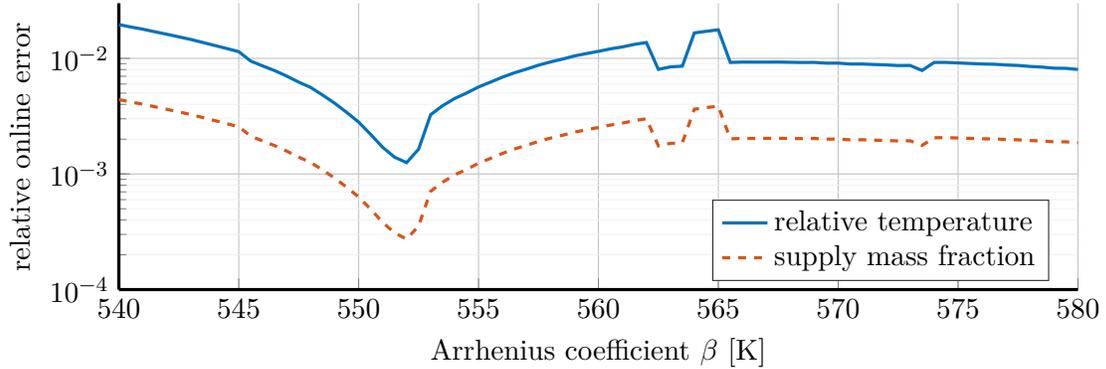

\begin{figure}
	\centering
	\input{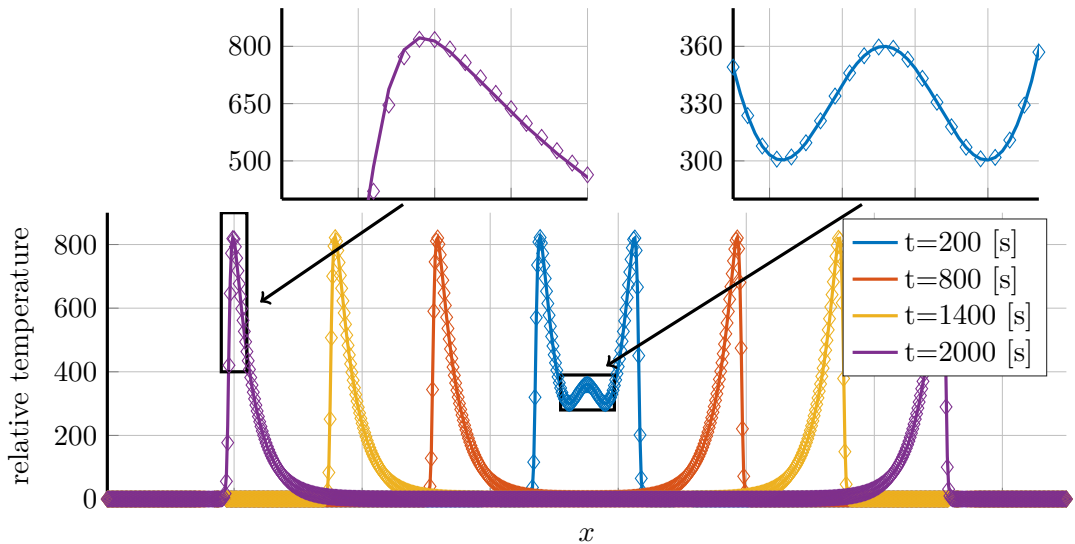}
	\caption{Solution of the temperature for the FOM (solid lines) and sPOD-sDEIM ROM (diamonds) at different time instances with Arrhenius coefficient $\beta=\si{540\kelvin}$}
	\label{fig:CombustionWavesFOMvsROM}
\end{figure}

Similarly as before, we perform a parameter study for the sPOD-sDEIM ROM by changing the Arrhenius coefficient $\beta$. 
To this end, we construct a ROM with $28$ POD modes before the switching time $\switchTime = \si{100\second}$, using $14$ modes for the relative temperature and $14$ modes for the supply mass fraction. 
After the switch, we use $3$ transformed modes for each variable and each frame and $2$ additional POD modes per variable, giving a total of $18$ DOF, including the path variables for each frame.
The modes are determined from snapshot data obtained from the simulation of the FOM for the training parameters $\beta\in\{\si{540\kelvin}, \si{560\kelvin}, \si{580\kelvin}\}$. 
The relative errors are depicted in \Cref{fig:parameterSampling_Gaussian}. 
Clearly, our ROM can produce accurate results throughout the parameter domain with an average speedup factor higher than $25$, see \Cref{tab:parameterStudy} for further details.
\begin{table}
	\centering
	\caption{Parameter study for the sPOD-sDEIM ROM over the Arrhenius coefficients with $81$ samples from $\beta=\si{540\kelvin}$ to $\beta=\si{580\kelvin}$}
	\label{tab:parameterStudy}
	\begin{tabular}{lrrr}
		\toprule
		& \textbf{min} & \textbf{mean} & \textbf{max}\\\midrule
		\textbf{relative error} & \num{1.25411e-03} & \num{9.18012e-03} & \num{ 1.95968e-02}\\
		\textbf{speedup} & \num{16.15} & \num{26.66} & \num{34.42}\\\bottomrule
	\end{tabular}
\end{table}
The worst-case approximation error throughout the parameter domain is achieved at $\beta=\si{540\kelvin}$. The corresponding solution of the FOM and the sPOD-sDEIM ROM for the temperature at different time instances is depicted in \Cref{fig:CombustionWavesFOMvsROM}, detailing that the wavefronts are captured accurately with a negligible offset.

\section{Conclusions}

In this paper, we introduce a hyper-reduction framework for the recently introduced nonlinear model reduction scheme presented in \cite{BlaSU20} and apply the method to a wildland fire model with nonlinear dynamics, featuring traveling combustion waves.
As approximation ansatz, we use a linear combination of transformed modes, cf.~\eqref{eqn:transformedMORansatz} and \eqref{eqn:transformedMORansatz:2}. Thus, the reduced state vector consists of the coefficients of the linear combination and the so-called path variables that parameterize the transformation operators.
The hyper-reduction framework addresses two sources of nonlinearity in the reduced-order model: First, we consider those nonlinearities originating from the nonlinear approximation ansatz.
Primarily, we propose to sample the path-dependent coefficient matrices in the offline phase and replace them with interpolation-based approximants in the online stage.
Second, we consider the nonlinearity originating from the FOM nonlinearity. We propose an extension of the discrete empirical interpolation method (DEIM) by approximating the nonlinearity with a linear combination of transformed modes.
All in all, this new hyper-reduction framework allows obtaining parameter-dependent ROMs based on transformed modes while achieving an efficient offline/online decomposition.

We apply the new method to two wildland fire test cases with different initial conditions, using a shift operator as a transformation to capture the traveling combustion waves.
The first considered initial condition consists of already separated combustion waves, such that the dynamic is mainly characterized by the propagation of the waves.
In this setting, the proposed model reduction framework based on transformed modes yields very low-dimensional ROMs, which clearly outperform ROMs constructed by a classical POD-DEIM approach in terms of accuracy and speedup.
Furthermore, the ROMs are parameter-dependent and prove to yield accurate predictions in the considered range of parameter values.
The second considered test case also includes the formation of the combustion waves, thus combining traveling and non-transported effects. Consequently, this is a challenging problem, not only for standard POD-DEIM ROMs, but also for a ROM based only on transported modes.
For this setting, we propose to take a low-dimensional POD-DEIM ROM for the initial time period where the combustion waves develop, and a ROM based on transformed modes and POD modes in the second part of the time interval.
Again, the ROM obtained by this approach is significantly faster and more accurate than classical POD-DEIM ROMs. Furthermore, the ROM yields accurate predictions within the considered parameter range.

Despite the convincing observed performance of our approach, there are still some shortcomings and difficulties which need to be addressed in the future.
First, we use the shift operator for the presented numerical experiments on a one-dimensional domain, but the extension to higher space dimensions is non-trivial.
To be able to cope with complex transports in higher dimensions, one promising direction is to use more involved coordinate transforms, see for instance \cite{KraSR21,RimPM20,Tad20} for some contributions in this direction.
Another difficulty that we observed in the numerical experiments is that the quality of the ROM obtained by our approach relies on a sufficiently fine spatial discretization of the FOM.
We assume that the reason for this is that we need to discretize the derivative of the shift operator for constructing the ROM and the corresponding discretization error scales with the spatial discretization error of the FOM, see also \cite[sec.~7.2]{BlaSU20}.
Furthermore, in contrast to the classical POD-DEIM approach, our method requires choosing more parameters. Their influence on the error is less evident than for POD-DEIM.

Besides the shortcomings mentioned in the last paragraph, an interesting future research direction is to ensure stability within the sPOD-sDEIM ROM based on transformed modes, for instance, by using a port-Hamiltonian formulation of the governing equations.
First investigations suggest that the considered wildland fire model allows for such an energy-based formulation, which could be exploited for structure-preserving model reduction schemes.
Furthermore, applying the proposed method to more involved wildland fire models is another interesting and challenging task for the future.
For instance, we could consider wildland fires formulated on two-dimensional domains or take into account the fluid dynamics in the atmosphere.

\subsection*{Acknowledgments}
The work of F.~Black is supported by the Deutsche Forschungsgemeinschaft (DFG, German Research Foundation) Collaborative Research Center 1029 \emph{Substantial efficiency increase in gas turbines through direct use of coupled unsteady combustion and flow dynamics}, project number 200291049. 
P.~Schulze acknowledges funding by the Deutsche Forschungsgemeinschaft (DFG, German Research Foundation) Collaborative Research Center Transregio 154 \emph{Mathematical Modelling, Simulation and Optimization Using the Example of Gas Networks }, project number 239904186. 
B. Unger acknowledges funding from the DFG under Germany's Excellence Strategy -- EXC 2075 -- 390740016 and is thankful for support by the Stuttgart Center for Simulation Science (SimTech). 

\bibliographystyle{plain}
\bibliography{literature}

\end{document}

%% file: def.tex
\DeclareMathOperator{\spann}{span}
\DeclareMathOperator{\diag}{diag}

\newcommand{\MOR}[1]{\widetilde{#1}}
\newcommand{\ROM}[1]{\widehat{#1}}

\newcommand{\N}{\ensuremath\mathbb{N}}

\newcommand{\R}{\ensuremath\mathbb{R}}

\newcommand{\calF}{\ensuremath\mathcal{F}}

\newcommand{\calT}{\ensuremath\mathcal{T}}

\newcommand{\calZ}{\ensuremath\mathcal{Z}}
\newcommand{\ansatzSpace}{\calZ}
\newcommand{\MORspace}{\MOR{\ansatzSpace}}

\newcommand{\arrhenius}{\beta}
\newcommand{\blkdiag}{\mathrm{blkdiag}}
\newcommand{\fuelRate}{\gamma_{\mathrm{S}}}
\newcommand{\heatTransfer}{\gamma}
\newcommand{\length}{\ell}
\newcommand{\Ltwo}{L^2}
\newcommand{\tempRise}{\alpha}

\newcommand{\state}{z}
\newcommand{\stateDim}{n}
\newcommand{\timeInt}{\ensuremath\mathbb{T}}
\newcommand{\finalTime}{t_\mathrm{f}}
\newcommand{\switchTime}{t_{\mathrm{switch}}}
\newcommand{\nSpacePoints}{n_x}
\newcommand{\param}{\mu}
\newcommand{\paramSet}{\mathbb{M}}

\newcommand{\ROMstate}{\ROM{\state}}
\newcommand{\ROMapprox}{\widetilde{\state}}
\newcommand{\ROMdim}{r}
\newcommand{\Ared}{\ROM{A}}
\newcommand{\fred}{\ROM{F}}

\newcommand{\DEIM}{{\mathrm{DEIM}}}
\newcommand{\sDEIM}{{\mathrm{sDEIM}}}
\newcommand{\DEIMcoeff}{b}
\newcommand{\DEIMdim}{m}
\newcommand{\Path}{p}

%% file: ROM-separatedWaves-modes-var1.tex
%
\definecolor{mycolor1}{rgb}{0.00000,0.44700,0.74100}%
\definecolor{mycolor2}{rgb}{0.85000,0.32500,0.09800}%
\begin{tikzpicture}

\begin{axis}[%
width=.8\linewidth,
height=1.5in,
at={(1.011in,0.642in)},
scale only axis,
grid=both,
grid style={line width=.1pt, draw=gray!10},
major grid style={line width=.2pt,draw=gray!50},
axis lines*=left,
axis line style={line width=\lineWidth},
xmin=0,
xmax=3000,
xmajorticks=false,
xlabel = {$x$},
ymin=-0.15,
ymax=0.15,
axis background/.style={fill=white},
legend style={legend cell align=left, align=left, draw=white!15!black,font=\small},
legend pos={north east},
reverse legend,
]

\addplot [color=mycolor2, line width=\lineWidth]
  table[row sep=crcr]{%
1	-0\\
1055	-6.61028279864695e-05\\
1059	-0.000206629015792714\\
1062	-0.00048575776099824\\
1064	-0.000858826395415235\\
1065	-0.00114195294008823\\
1066	-0.0015184168637461\\
1067	-0.00201898848536075\\
1068	-0.0026845818692891\\
1069	-0.00356959937607826\\
1070	-0.00474637761681151\\
1071	-0.00631110069753049\\
1072	-0.00839166067544284\\
1073	-0.0111581147039033\\
1074	-0.0148366234861896\\
1075	-0.019727996052552\\
1076	-0.0262296339074055\\
1077	-0.0348351634102073\\
1078	-0.0459983978739729\\
1079	-0.0596890366032312\\
1080	-0.0747818774502775\\
1081	-0.0889790337541854\\
1082	-0.0997077177034953\\
1083	-0.10530259690313\\
1084	-0.105486110669062\\
1085	-0.101051584161723\\
1086	-0.093274927703078\\
1087	-0.0834562027498578\\
1088	-0.0726848450967736\\
1089	-0.0617700363145559\\
1090	-0.0512578303446389\\
1091	-0.0414813774632421\\
1092	-0.0326160764825545\\
1093	-0.0247277940875392\\
1094	-0.0178107349215679\\
1095	-0.0118153240682659\\
1096	-0.00666778809954849\\
1097	-0.00228336489544745\\
1098	0.0014251556885938\\
1099	0.0045421815152622\\
1100	0.00714633107236295\\
1101	0.00930894596740472\\
1102	0.011093578287273\\
1103	0.0125560730530196\\
1104	0.0137450053020984\\
1105	0.0147022952132829\\
1106	0.0154639050251717\\
1107	0.016060539636328\\
1108	0.0165183231897572\\
1109	0.0168594175856924\\
1111	0.0172636861557294\\
1113	0.0173911358319856\\
1115	0.0173256066987051\\
1118	0.0169908708571711\\
1122	0.0162936503998026\\
1129	0.0148131354440011\\
1140	0.0124827695972272\\
1148	0.0109621655938099\\
1156	0.00961208707758487\\
1164	0.00842274338629068\\
1173	0.00725710656797673\\
1183	0.00614887014990018\\
1193	0.00520930185393809\\
1204	0.00434041071139291\\
1216	0.00355661285857423\\
1230	0.00281848574786636\\
1246	0.00215807640051935\\
1264	0.00158838617880974\\
1290	0.000951500347127876\\
1314	0.00031207868232741\\
1326	-0.000180712404926453\\
1335	-0.000715515350748319\\
1342	-0.00127715575763432\\
1349	-0.00201363008409317\\
1355	-0.00281985160472686\\
1360	-0.00364172393210538\\
1365	-0.00462685413413055\\
1370	-0.0057963036797446\\
1375	-0.00718012013658154\\
1379	-0.00845751558063057\\
1383	-0.00990241687168236\\
1387	-0.0115290393400755\\
1391	-0.0133513758573827\\
1394	-0.0148533355254585\\
1398	-0.0170501598745432\\
1402	-0.0194778178456545\\
1405	-0.0214580861952527\\
1408	-0.023575019427426\\
1411	-0.0258383153309296\\
1414	-0.0282440856553876\\
1417	-0.030801052231709\\
1420	-0.0335034083309438\\
1423	-0.0363586459316139\\
1426	-0.0393588127553812\\
1429	-0.0425099685871828\\
1432	-0.0458018573253867\\
1436	-0.0504125652187213\\
1440	-0.0552632038684351\\
1444	-0.0603397012114328\\
1448	-0.0656250082120096\\
1452	-0.0710993058173699\\
1456	-0.0767402991441486\\
1461	-0.0839910756894824\\
1466	-0.09140812753094\\
1474	-0.103501871753906\\
1487	-0.123246904076041\\
1494	-0.133683476653005\\
1499	-0.140994824722384\\
1500	-0.142436958924009\\
1501	-0.133311196981595\\
1502	-0.125992162690181\\
1503	-0.118418952799402\\
1504	-0.110972556311026\\
1505	-0.103699132244401\\
1506	-0.0967701069871509\\
1507	-0.0891575676050707\\
1508	-0.0826824673549709\\
1509	-0.0760907029857663\\
1510	-0.0697026208313218\\
1511	-0.0635733052254182\\
1512	-0.0577061734575182\\
1513	-0.0514850015647426\\
1515	-0.0408794518102695\\
1516	-0.0357964599702427\\
1517	-0.0309677370864847\\
1518	-0.0263247146840513\\
1519	-0.0215748570913092\\
1521	-0.0134925011921041\\
1522	-0.00966517558936175\\
1524	-0.0026229869376948\\
1525	0.00079616962966611\\
1527	0.00661119907499597\\
1528	0.00932770653162152\\
1530	0.0142692839822303\\
1531	0.0165874657027416\\
1533	0.0205180376656244\\
1534	0.0223141080241476\\
1536	0.0255530820995773\\
1537	0.0270104374180846\\
1539	0.0294578647667549\\
1540	0.0305470462371886\\
1542	0.0324725286977809\\
1543	0.0332849919027467\\
1545	0.0346167290185804\\
1546	0.0351828983589257\\
1549	0.0364869239879226\\
1551	0.03701497429347\\
1553	0.0373201602401423\\
1554	0.0374704424971242\\
1556	0.0375040050139432\\
1558	0.0374087568957293\\
1561	0.0370300389031399\\
1564	0.036403022378181\\
1568	0.0352898385531262\\
1572	0.033953229470626\\
1579	0.031260284551081\\
1599	0.023282682451736\\
1605	0.0211139531493245\\
1611	0.0191046489558175\\
1617	0.0172605658522116\\
1623	0.0155793627518506\\
1630	0.0138147060411029\\
1636	0.0124573288503598\\
1643	0.0110386224314425\\
1651	0.00961664435317289\\
1659	0.00837903758701941\\
1667	0.00730158151372962\\
1676	0.00625658884564473\\
1686	0.00527096275800432\\
1697	0.00436591925199536\\
1709	0.00355501145395465\\
1722	0.00284503358670918\\
1737	0.00219894272549936\\
1754	0.00164087214670872\\
1773	0.00118082756671356\\
1797	0.000777425246269559\\
1827	0.000458891648577264\\
1866	0.00022922582411411\\
1922	8.28758702482446e-05\\
2025	1.1705268661899e-05\\
2457	0\\
3000	0\\
};
\addlegendentry{second mode}

\addplot [color=mycolor1, line width=\lineWidth]
  table[row sep=crcr]{%
1	-0\\
1059	7.9590296081733e-05\\
1063	0.000248789122906601\\
1065	0.000439862586517847\\
1067	0.000777683095293469\\
1069	0.0013749543945778\\
1070	0.001828231178024\\
1071	0.00243093825429241\\
1072	0.0032323378077308\\
1073	0.00429793215380414\\
1074	0.00571481968745502\\
1075	0.00759883282535156\\
1076	0.0101040149115761\\
1077	0.0134339645528598\\
1078	0.0178473077962735\\
1079	0.0236316963982972\\
1080	0.0310151130347549\\
1081	0.0400329110520943\\
1082	0.0504421240389092\\
1083	0.0617596415545449\\
1084	0.0733960114139336\\
1085	0.0847946705621325\\
1086	0.0955170469919722\\
1087	0.105268518025696\\
1088	0.113886544533216\\
1089	0.121312797245082\\
1090	0.127562947853676\\
1091	0.132700460837896\\
1092	0.136816302797342\\
1093	0.140014367016192\\
1094	0.142401622558737\\
1095	0.144081856196863\\
1096	0.145152000331564\\
1097	0.145700243531337\\
1098	0.145805321230455\\
1099	0.145536553167858\\
1100	0.14495432569629\\
1101	0.144110814691885\\
1102	0.143050814948765\\
1103	0.141812591071812\\
1104	0.140428698585765\\
1105	0.138926746576089\\
1107	0.135658433102435\\
1109	0.132153919007578\\
1112	0.12667213880286\\
1118	0.115613636729904\\
1121	0.110247144576988\\
1124	0.105056853014048\\
1127	0.100064429273516\\
1130	0.095279221338842\\
1133	0.0907031931551501\\
1136	0.0863339075012846\\
1139	0.0821663413921669\\
1142	0.0781939985540703\\
1145	0.0744095952677526\\
1148	0.0708054850097142\\
1151	0.0673739217463662\\
1154	0.0641072226576398\\
1157	0.0609978676629908\\
1160	0.0580385588887111\\
1163	0.0552222546134544\\
1166	0.0525421868328522\\
1170	0.0491695146183702\\
1174	0.0460129385896835\\
1178	0.0430587255764294\\
1182	0.040293979186572\\
1186	0.0377066013193144\\
1190	0.0352852507016905\\
1194	0.0330193006811896\\
1198	0.0308987975627133\\
1203	0.0284385602972179\\
1208	0.0261741524100216\\
1213	0.0240900019089167\\
1218	0.0221717712088321\\
1223	0.0204062604007049\\
1229	0.0184721966843426\\
1235	0.0167214236021209\\
1241	0.0151365769097538\\
1248	0.0134764119052306\\
1255	0.0119983285771923\\
1262	0.0106823589835585\\
1270	0.00935418434301027\\
1278	0.00819114707928748\\
1287	0.00705465865530641\\
1297	0.0059758540310213\\
1308	0.00497871417792339\\
1320	0.00407970786909573\\
1333	0.00328805671279042\\
1347	0.00260649875554009\\
1363	0.00199897589072862\\
1382	0.00145912415791827\\
1404	0.00101454748164542\\
1431	0.000652174261631444\\
1465	0.000379808073830645\\
1515	0.000172819050476392\\
1582	5.164010099179e-05\\
1726	3.71842997992644e-06\\
3000	0\\
};
\addlegendentry{first mode}

\end{axis}

\end{tikzpicture}%

%% file: ROM-separatedWaves-modes-var2.tex
%
\definecolor{mycolor1}{rgb}{0.00000,0.44700,0.74100}%
\definecolor{mycolor2}{rgb}{0.85000,0.32500,0.09800}%
\begin{tikzpicture}

\begin{axis}[%
width=.8\linewidth,
height=1.5in,
at={(4.436in,0.642in)},
scale only axis,
grid=both,
grid style={line width=.1pt, draw=gray!10},
major grid style={line width=.2pt,draw=gray!50},
axis lines*=left,
axis line style={line width=\lineWidth},
xmin=0,
xmax=3000,
xmajorticks=false,
xlabel={$x$},
ymin=-0.15,
ymax=0.15,
axis background/.style={fill=white},
legend style={legend cell align=left, align=left, draw=white!15!black,font=\small},
legend pos={north east},
reverse legend,
]

\addplot [color=mycolor2, line width=\lineWidth]
  table[row sep=crcr]{%
1	-0.00203865806543035\\
1076	-0.00199594196919861\\
1077	-0.00127346557155761\\
1078	0.00310706088021107\\
1079	0.0183513204888186\\
1080	0.0532024894087044\\
1081	0.109724129676579\\
1082	0.178520438319083\\
1083	0.244138222675701\\
1084	0.29417122562063\\
1085	0.323461949958528\\
1086	0.332918158148914\\
1087	0.326566106115024\\
1088	0.309206818305029\\
1089	0.285162984839644\\
1090	0.257822785259577\\
1091	0.229609291541692\\
1092	0.202126092905019\\
1093	0.17634747932334\\
1094	0.152795163471183\\
1095	0.131681931398816\\
1096	0.113020154921742\\
1097	0.0966998712679015\\
1098	0.0825428755115354\\
1099	0.0703389151708507\\
1100	0.0598690178608194\\
1101	0.0509198315271533\\
1102	0.0432918453766433\\
1103	0.036803550865443\\
1104	0.0312929882770732\\
1105	0.0266176738064132\\
1106	0.0226535789283844\\
1107	0.0192936063149318\\
1108	0.0164458488466153\\
1109	0.0140318105486585\\
1110	0.01198469584142\\
1111	0.0102478256071663\\
1112	0.00877320761856026\\
1113	0.00752026941154327\\
1114	0.00645474994689721\\
1115	0.00554773977955847\\
1117	0.00411553815956722\\
1119	0.00307095963535176\\
1121	0.00230535999435233\\
1123	0.00174131933863464\\
1126	0.00115661092922892\\
1130	0.000685311799315969\\
1135	0.000369186893294682\\
1143	0.000148769592669851\\
1158	3.46633214576286e-05\\
1209	8.66216396389063e-07\\
1397	0.000101651590739493\\
1416	0.00033772011011024\\
1432	0.000768518127642892\\
1461	0.00187111858758726\\
1506	0.00348706996828696\\
1549	0.00480378759175437\\
1592	0.00589474326943673\\
1635	0.00675897927521873\\
1678	0.00739573929922699\\
1721	0.00780446689532255\\
1764	0.00798480093180842\\
1807	0.00793658474458425\\
1850	0.00765986207397873\\
1893	0.00715487167235551\\
1936	0.0064220583790302\\
1979	0.00546206732542487\\
2022	0.00427581153462597\\
2066	0.00283096850489528\\
2143	0.000214864027384465\\
2159	3.93056993743812e-05\\
2193	9.24330379348248e-07\\
3000	0\\
};
\addlegendentry{second mode}

\addplot [color=mycolor1, line width=\lineWidth]
  table[row sep=crcr]{%
1	0.0303476178105484\\
1079	0.030250627166879\\
1080	0.0300042596390995\\
1081	0.0294427652988816\\
1082	0.0284441425314981\\
1083	0.0269767479517213\\
1084	0.0251086433027012\\
1085	0.022968618892719\\
1088	0.0162293368425708\\
1089	0.0141823617041155\\
1090	0.0123140600535407\\
1091	0.0106381199148018\\
1092	0.00915452831623043\\
1093	0.00785454639071759\\
1094	0.0067244051497255\\
1095	0.00574789683378185\\
1096	0.00490810236760808\\
1098	0.00357351668708361\\
1100	0.00260238089276754\\
1102	0.001899192113342\\
1104	0.00139087709112573\\
1107	0.000879325403730036\\
1111	0.000486621750496852\\
1116	0.000240864937495644\\
1125	7.59775916776562e-05\\
1146	9.29062707655248e-06\\
1318	9.03582531464053e-07\\
3000	0\\
};
\addlegendentry{first mode}
\end{axis}

\end{tikzpicture}%

%% file: separatedWaves-POD-vs-SPOD.tex
%
\definecolor{mycolor1}{rgb}{0.00000,0.44700,0.74100}%
\definecolor{mycolor2}{rgb}{0.85000,0.32500,0.09800}%
\definecolor{mycolor3}{rgb}{0.92900,0.69400,0.12500}%
\begin{tikzpicture}

\begin{axis}[%
width=.85\linewidth,
height=1.7in,
at={(0.758in,0.481in)},
scale only axis,
grid=both,
grid style={line width=.1pt, draw=gray!10},
major grid style={line width=.2pt,draw=gray!50},
axis lines*=left,
axis line style={line width=\lineWidth},
xmin=0,
xmax=0.4,
xlabel style={font=\color{white!15!black}},
x tick label style = {/pgf/number format/.cd, fixed},
xlabel={computation time [s]},
ymode=log,
ymin=0.0001,
ymax=10,
yminorticks=true,
ylabel style={font=\color{white!15!black}},
ylabel={relative online error},
axis background/.style={fill=white},
title style={font=\bfseries},
legend style={legend cell align=left, align=left, draw=white!15!black,font=\small}
]
\addplot [color=mycolor1, line width=1.5, only marks, mark size=3pt ,mark=+, mark options={solid, mycolor1}]
  table[row sep=crcr]{%
0.014010541	0.00617727113246551\\ 
0.02292975	0.000442839144601316\\ 
0.044229541	0.000243302627710679\\ 
0.05315275	0.000258429160162966\\ 
0.073627958	0.000190398075805284\\ 
0.095732541	0.00028293562122469\\ 
0.121406042	0.000257030841564541\\ 
0.169113083	0.000376038791697564\\ 
};
\node[font=\tiny, pin= above right:{\small 4+2}] at (axis cs:0.014010541,0.00617727113246551) {};
\node[font=\tiny, pin= above right:{\small 8+2}] at (axis cs:0.02292975	, 0.000442839144601316) {};
\node[font=\tiny, pin= above right:{\small 20+2}] at (axis cs:0.073627958, 0.000190398075805284) {};
\node[font=\tiny, pin= above right:{\small 32+2}] at (axis cs:0.169113083,0.000376038791697564) {};
\addlegendentry{sPOD-sDEIM}

\addplot [color=mycolor2, line width=1.5, only marks, mark size=2pt ,mark=*, mark options={solid, mycolor2} ]
  table[row sep=crcr]{%
0.0144217916666667	1.05543511610239\\ 
0.0175923473333333	0.881027957035218\\ 
0.024540514	0.697757665433012\\ 
0.0309287086666667	0.525721583110737\\ 
0.192945722	0.377497492370512\\ 
0.24386493	0.259251929886418\\ 
0.283281430666667	0.173727295082382\\ 
0.337009833333333	0.11402056942729\\ 
};
\node[font=\tiny, pin= above right:{\small 60}] at (axis cs:0.0144217916666667,1.05543511610239) {};
\node[font=\tiny, pin= above right:{\small 120}] at (axis cs:0.0309287086666667, 0.525721583110737) {};
\node[font=\tiny, pin= above right:{\small 140}] at (axis cs:0.192945722, 0.377497492370512) {};
\node[font=\tiny, pin= below right:{\small 200}] at (axis cs:0.337009833333333, 0.11402056942729) {};
\addlegendentry{POD-DEIM}

\end{axis}

\end{tikzpicture}%


%% file: separatedWavesParameterSamplingN3000r3-3hyper.tex
%
%
\definecolor{mycolor1}{rgb}{0.00000,0.44700,0.74100}%
\definecolor{mycolor2}{rgb}{0.85000,0.32500,0.09800}%
\definecolor{mycolor3}{rgb}{0.92900,0.69400,0.12500}%
\begin{tikzpicture}

\begin{axis}[%
width=.85\linewidth,
height=1.5in,
at={(1.011in,0.642in)},
scale only axis,
grid=both,
grid style={line width=.1pt, draw=gray!10},
major grid style={line width=.2pt,draw=gray!50},
axis lines*=left,
axis line style={line width=\lineWidth},
xmin=540,
xmax=580,
xlabel={Arrhenius coefficient $\arrhenius$ [\si{\kelvin}]},
ymin=1e-4,
ymax=1e-2,
ymode=log,
ylabel={relative online error},
axis background/.style={fill=white},
legend style={legend cell align=left, align=left, draw=white!15!black,font=\small},
legend pos={south east},
]

\addplot [color=mycolor1,line width=\lineWidth]
  table[row sep=crcr]{%
540	0.00612006901883122\\
541	0.00575130886299666\\
542	0.00539069727342518\\
543	0.00503575092834216\\
544.5	0.00451973677779028\\
545.5	0.00416795347574796\\
546	0.00399449776020974\\
547	0.00364346311965844\\
547.5	0.00346528154241057\\
548	0.00328854056883426\\
548.5	0.00311108173436714\\
549	0.00293029069983725\\
549.5	0.00274885616628725\\
550	0.00256407752595878\\
550.5	0.00239749711775199\\
551	0.00221050941256195\\
551.5	0.00202058775551345\\
552	0.00182555243465547\\
552.5	0.00161360455981567\\
553	0.00141916756300086\\
553.5	0.00122340360940256\\
554	0.00102855898975232\\
554.5	0.000843098947914242\\
555	0.000676074849475488\\
555.5	0.00055190242318908\\
556	0.000513393604635335\\
556.5	0.000584826960313096\\
557	0.000744116738092248\\
557.5	0.000954350678785297\\
558	0.00119329852352563\\
558.5	0.0014550839698209\\
559	0.00173415763737191\\
559.5	0.00202899732791332\\
560	0.00232175201588754\\
560.5	0.00264209255258334\\
561	0.00296169714045279\\
561.5	0.00330423065203625\\
562	0.00364823310096103\\
562.5	0.00401154466397876\\
563	0.00435818943454038\\
563.5	0.00470447176113958\\
564	0.00505145127571594\\
564.5	0.00537527093100269\\
565	0.00569125666723633\\
565.5	0.005983828306622\\
566	0.00624487573104274\\
566.5	0.0064860716448524\\
567	0.00671149710138547\\
567.5	0.00691207828184272\\
568	0.00708532509803808\\
568.5	0.00723573508302477\\
569	0.00735993795787367\\
569.5	0.00747574167784576\\
570	0.00756486570426149\\
570.5	0.00763943636335207\\
571.5	0.00773990328816062\\
572.5	0.0078068012772282\\
573.5	0.00784078574533782\\
574.5	0.00787934369997437\\
575.5	0.00792770765549399\\
576.5	0.00800795450435875\\
577	0.00806660396280752\\
577.5	0.00814658102739242\\
578	0.00823965230912256\\
578.5	0.00835362539673951\\
579	0.00850112540824264\\
579.5	0.00866919805311005\\
580	0.00888574607661374\\
};
\addlegendentry{relative temperature}

\addplot [color=mycolor2,dashed,line width=\lineWidth]
  table[row sep=crcr]{%
540	0.00135427006184932\\
541	0.00127237282189752\\
542	0.00119249636665854\\
543	0.00111408592415355\\
544.5	0.00100060769569631\\
545.5	0.000923164573199\\
546	0.000885099556085458\\
547	0.000807978180915739\\
547.5	0.00076873991063812\\
548	0.000729918715376931\\
548.5	0.00069096217731551\\
549	0.000651220961302404\\
549.5	0.000611326281128135\\
550	0.000570674042009035\\
550.5	0.000534485205708989\\
551	0.000493345925983999\\
551.5	0.000451597049723161\\
552	0.000408681296863528\\
552.5	0.000361795422905248\\
553	0.00031922355089309\\
553.5	0.000276369532590064\\
554	0.000233943246963788\\
554.5	0.00019400801873586\\
555	0.000158693309414373\\
555.5	0.000133228691131835\\
556	0.000125969924747171\\
556.5	0.000141357645664499\\
557	0.000175385979317608\\
557.5	0.00022099264195878\\
558	0.000273534686870933\\
558.5	0.000331504606543174\\
559	0.000393603426284655\\
559.5	0.000459437441501873\\
560	0.00052488903464083\\
560.5	0.000596737609833722\\
561	0.000668446385960688\\
561.5	0.000745457338283697\\
562	0.000822804133787026\\
562.5	0.000904626941271803\\
563	0.000982622699875204\\
563.5	0.00106059664954058\\
564	0.00113876835538355\\
564.5	0.00121168665287583\\
565	0.00128289215431788\\
565.5	0.00134882306778995\\
566	0.00140760111803622\\
566.5	0.00146197313276061\\
567	0.00151286148695925\\
567.5	0.00155830785204974\\
568	0.00159764378590709\\
568.5	0.00163197513500558\\
569.5	0.00168733605455498\\
570	0.0017083682709052\\
570.5	0.00172635607832183\\
571.5	0.00175201036021366\\
572.5	0.00177072886239311\\
575.5	0.00181679065670052\\
576.5	0.00184337465614627\\
577	0.00186110026889053\\
577.5	0.00188399763416246\\
578	0.00191005321832946\\
578.5	0.0019409377449421\\
579	0.00197955264269211\\
579.5	0.0020226886966473\\
580	0.00207721876720349\\
};
\addlegendentry{supply mass fraction}

\end{axis}
\end{tikzpicture}%

%% file: varyNumbersPODModes_T_offline.tex
%
%
\definecolor{mycolor1}{rgb}{0.00000,0.44700,0.74100}%
\definecolor{mycolor2}{rgb}{0.85000,0.32500,0.09800}%
\definecolor{mycolor3}{rgb}{0.92900,0.69400,0.12500}%
\definecolor{mycolor4}{rgb}{0.49400,0.18400,0.55600}%
\definecolor{mycolor5}{rgb}{0.46600,0.67400,0.18800}%

\begin{tikzpicture}

\begin{axis}[%
width=.75\linewidth,
height=1.5in,
at={(1.081in,0.864in)},
scale only axis,
xmin=1,
xmax=41,
xlabel style={font=\color{white!15!black}},
xlabel={\# POD modes in 1st time interval},
ymode=log,
ymin=0.0001,
ymax=0.1,
grid=both,
grid style={line width=.1pt, draw=gray!10},
major grid style={line width=.2pt,draw=gray!50},
yminorticks=true,
axis line style={line width=\lineWidth},
ylabel style={font=\color{white!15!black}},
ylabel={relative offline error},
axis background/.style={fill=white},
axis x line*=bottom,
axis y line*=left,
legend style={legend cell align=left, align=left, draw=white!15!black,column sep=1ex, at={(2.1,1.35)},font=\small},
legend columns=2,
]
\addplot [color=mycolor1, line width=1.5pt, only marks, mark size=3pt, mark=+, mark options={solid, mycolor1}]
  table[row sep=crcr]{%
4	0.0389904277657001\\
8	0.0125409207892634\\
12	0.00561243519181379\\
16	0.00405200339444851\\
20	0.00377661800288958\\
24	0.00373000632785823\\
28	0.00372170576180886\\
32	0.00372015307295421\\
36	0.00371985318607653\\
40	0.00371979426182288\\
};
\addlegendentry{2 POD modes in 2nd time interval}


\addplot [color=mycolor3, line width=1.5pt, only marks, mark size=2pt, mark=*, mark options={solid, mycolor2}]
  table[row sep=crcr]{%
4	0.0388139660850303\\
8	0.011981035019571\\
12	0.00421544003396049\\
16	0.00163988944939257\\
20	0.000730307768927908\\
24	0.000428314257840286\\
28	0.000348711449608571\\
32	0.000331729986310324\\
36	0.000328349849353555\\
40	0.000327681625815723\\
};
\addlegendentry{6 POD modes in 2nd time interval}


\addplot [color=mycolor5, line width=1.5pt, only marks, mark size=3pt, mark=star, mark options={solid, mycolor3}]
  table[row sep=crcr]{%
4	0.0388132261749929\\
8	0.0119786377773709\\
12	0.00420862181018766\\
16	0.00162228242531919\\
20	0.000689863971103648\\
24	0.000354987274053186\\
28	0.000253303249798632\\
32	0.00022936356530165\\
36	0.000224447064674366\\
40	0.00022346836278917\\
};
\addlegendentry{10 POD modes in 2nd time interval}

\end{axis}
x
\end{tikzpicture}%

%% file: varyNumbersPODModes_S_offline.tex
%
%
\definecolor{mycolor1}{rgb}{0.00000,0.44700,0.74100}%
\definecolor{mycolor2}{rgb}{0.85000,0.32500,0.09800}%
\definecolor{mycolor3}{rgb}{0.92900,0.69400,0.12500}%
\definecolor{mycolor4}{rgb}{0.49400,0.18400,0.55600}%
\definecolor{mycolor5}{rgb}{0.46600,0.67400,0.18800}%

\begin{tikzpicture}

\begin{axis}[%
width=.75\linewidth,
height=1.5in,
at={(1.081in,0.864in)},
scale only axis,
xmin=1,
xmax=41,
grid=both,
grid style={line width=.1pt, draw=gray!10},
major grid style={line width=.2pt,draw=gray!50},
xlabel style={font=\color{white!15!black}},
xlabel={\# POD modes in 1st time interval},
ymode=log,
ymin=1e-05,
ymax=0.0110407752041685,
yminorticks=true,
axis line style={line width=\lineWidth},
ylabel style={font=\color{white!15!black}},
ylabel={relative offline error},
axis background/.style={fill=white},
axis x line*=bottom,
axis y line*=left,
legend style={legend cell align=left, align=left, draw=white!15!black,column sep=1ex,font=\small},
legend columns=-1,
]
\addplot [color=mycolor1, line width=1.5pt, only marks, mark size=3pt, mark=+, mark options={solid, mycolor1}]
  table[row sep=crcr]{%
4	0.0110407752041685\\
8	0.00334364641071162\\
12	0.00126441062069685\\
16	0.000532294682825792\\
20	0.000239687610220734\\
24	0.000114870163154641\\
28	6.11816368803315e-05\\
32	4.01089690173876e-05\\
36	3.33244419111318e-05\\
40	3.15039721126291e-05\\
};


\addplot [color=mycolor3, line width=1.5pt, only marks, mark size=2pt, mark=*, mark options={solid, mycolor2}]
  table[row sep=crcr]{%
4	0.0110407553553437\\
8	0.00334358086898374\\
12	0.0012642372899501\\
16	0.000531882822559619\\
20	0.000238771560414349\\
24	0.0001129462790558\\
28	5.74882619387086e-05\\
32	3.42116495512067e-05\\
36	2.59273214595188e-05\\
40	2.35416190660317e-05\\
};


\addplot [color=mycolor5, line width=1.5pt, only marks, mark size=3pt, mark=star, mark options={solid, mycolor3}]
  table[row sep=crcr]{%
4	0.0110407488843357\\
8	0.00334355950116688\\
12	0.00126418077651483\\
16	0.000531748481279683\\
20	0.000238472154500197\\
24	0.000112311942204795\\
28	5.62317585532947e-05\\
32	3.20553798024146e-05\\
36	2.3007312040473e-05\\
40	2.02809821634115e-05\\
};

\end{axis}
\end{tikzpicture}%

%% file: Gaussian_POD_vs_SPOD.tex
%
\definecolor{mycolor1}{rgb}{0.00000,0.44700,0.74100}%
\definecolor{mycolor2}{rgb}{0.85000,0.32500,0.09800}%
\definecolor{mycolor3}{rgb}{0.92900,0.69400,0.12500}%
\begin{tikzpicture}

\begin{axis}[%
width=.85\linewidth,
height=1.7in,
at={(0.758in,0.481in)},
scale only axis,
grid=both,
grid style={line width=.1pt, draw=gray!10},
major grid style={line width=.2pt,draw=gray!50},
axis lines*=left,
axis line style={line width=\lineWidth},
xmin=0,
xmax=1.5,
xlabel style={font=\color{white!15!black}},
xlabel={computation time [s]},
ymode=log,
ymin=0.001,
ymax=5,
yminorticks=true,
ylabel style={font=\color{white!15!black}},
ylabel={relative online error},
axis background/.style={fill=white},
title style={font=\bfseries},
legend style={legend cell align=left, align=left, draw=white!15!black,font=\small}
]
\addplot [color=mycolor1, line width=1.5, only marks, mark size=3pt ,mark=+, mark options={solid, mycolor1}]
  table[row sep=crcr]{%
0.050201708	0.217377128696895\\ 
0.081804459	0.210971646484724\\ 
0.130671208	0.203912559542678\\ 
0.054622042	0.217123980710477\\ 
0.089916709	0.210866159694033\\ 
0.133163292	0.203858025227596\\ 
0.19613425	0.00366434654022996\\ 
0.241465833	0.00489794344311111\\ 
0.25647	0.0033180240165287\\ 
0.210987834	0.00367173884174056\\ 
0.24729225	0.00529319801521602\\ 
0.277193709	0.0036085069830082\\ 
0.451444125	0.0089162067160263\\ 
0.611854417	0.0042114653287235\\ 
0.564206959	0.00229540170114418\\ 
0.482336791	0.00922477625651489\\ 
0.596952666	0.00466024873742848\\ 
0.560641083	0.00253509418798022\\ 
};
\node[font=\tiny, pin= below right:{\small 28--(4+4+2)}] at (axis cs:0.081804459, 0.210971646484724) {};
\node[font=\tiny, pin= above right:{\small 36--(12+6+2)}] at (axis cs:0.277193709, 0.0036085069830082) {};
\node[font=\tiny, pin= above right:{\small 36--(20+2+2)}] at (axis cs:0.482336791, 0.00922477625651489) {};
\node[font=\tiny, pin= above right:{\small 28--(20+4+2)}] at (axis cs:0.611854417, 0.0042114653287235) {};
\addlegendentry{sPOD-sDEIM}

\addplot [color=mycolor2, line width=1.5, only marks, mark size=2pt ,mark=*, mark options={solid, mycolor2} ]
  table[row sep=crcr]{%
0.0183932916666667	1.20221536265155\\ 
0.0196145693333333	1.14155125403256\\ 
0.024072444	1.06044580862066\\ 
0.032000125	0.96217757089737\\ 
0.228846514	0.85366271573802\\ 
0.263809360666667	0.741031825376094\\ 
0.309658902666667	0.629975611225554\\ 
0.352201028	0.52321037057425\\ 
0.402414291333333	0.426669867997251\\ 
0.472360222333333	0.339833475213526\\ 
0.552419444333333	0.266363857338503\\ 
0.653387528	0.205644784349842\\ 
0.793569236	0.156642729030763\\ 
1.077889	0.118790973147405\\ 
1.44713875	0.0896275463156704\\ 
1.695254806	0.0675576781057568\\ 
2.137017236	0.0506845301859442\\ 
2.61534266633333	0.0379585723105365\\ 
};
\node[font=\tiny, pin= above right:{\small 120}] at (axis cs:0.032000125, 0.96217757089737) {};
\node[font=\tiny, pin= above right:{\small 220}] at (axis cs:0.352201028, 0.52321037057425) {};
\node[font=\tiny, pin= above right:{\small 300}] at (axis cs:0.793569236, 0.156642729030763) {};
\node[font=\tiny, pin= below right:{\small 320}] at (axis cs:1.077889, 0.118790973147405) {};
\addlegendentry{POD-DEIM}

\end{axis}

\end{tikzpicture}%


%% file: ROM-Gaussian-parameterTesting.tex
%
\definecolor{mycolor1}{rgb}{0.00000,0.44700,0.74100}%
\definecolor{mycolor2}{rgb}{0.85000,0.32500,0.09800}%
\definecolor{mycolor3}{rgb}{0.92900,0.69400,0.12500}%
\begin{tikzpicture}

\begin{axis}[%
width=.85\linewidth,
height=1.5in,
at={(1.011in,0.642in)},
scale only axis,
grid=both,
grid style={line width=.1pt, draw=gray!10},
major grid style={line width=.2pt,draw=gray!50},
axis lines*=left,
axis line style={line width=\lineWidth},
xmin=540,
xmax=580,
xlabel={Arrhenius coefficient $\arrhenius$ [\si{\kelvin}]},
ymin=1e-4,
ymax=3e-2,
ymode=log,
ylabel={relative online error},
axis background/.style={fill=white},
legend style={legend cell align=left, align=left, draw=white!15!black},
legend pos={south east},
]

\addplot [color=mycolor1,line width=\lineWidth]
  table[row sep=crcr]{%
540	0.0195970314767195\\
540.5	0.0186799858752284\\
541	0.017882729441211\\
541.5	0.0169839039019649\\
542	0.0161664620474634\\
543	0.0145755669144513\\
544.5	0.0121647364712521\\
545	0.0114241880564196\\
545.5	0.00946692457505831\\
546	0.00859430706357953\\
546.5	0.00783487213187785\\
547	0.0070229402595106\\
547.5	0.00622935136668018\\
548	0.00560131834963787\\
548.5	0.00480945378779509\\
549	0.00409335602585667\\
549.5	0.00341072478433765\\
550	0.00280964859634324\\
550.5	0.00220108644361427\\
551	0.00170184380852062\\
551.5	0.00139557053697553\\
552	0.00125412838951879\\
552.5	0.00163866399700559\\
553	0.00325084586946565\\
553.5	0.00389324473864343\\
554	0.00449237651359869\\
554.5	0.00500948187616731\\
555	0.00565678060870574\\
555.5	0.00623020529358456\\
556	0.00688107170782981\\
556.5	0.00750549134721655\\
557	0.00807222387924505\\
557.5	0.0087279766573529\\
558	0.00931274108195372\\
558.5	0.00985147214494791\\
559	0.0104925581014229\\
559.5	0.0110000289446324\\
560	0.0115072403773873\\
560.5	0.0120845538209706\\
561	0.0125567671389016\\
561.5	0.0132168237388188\\
562	0.0137030868269106\\
562.5	0.00802661722263334\\
563	0.00844124069281372\\
563.5	0.00856108229926006\\
564	0.0166021252838537\\
564.5	0.0171452003955965\\
565	0.0176643916954388\\
565.5	0.00922234959147288\\
566	0.00928606333532696\\
567	0.00928000645376115\\
567.5	0.00929386894802143\\
568.5	0.00921093148569886\\
569	0.00922705438101979\\
569.5	0.00908933698745305\\
570	0.00909547506276113\\
570.5	0.00893573651084664\\
571	0.00893611833925024\\
571.5	0.00884468853999877\\
572	0.0087771015035446\\
572.5	0.00866244662263262\\
573	0.00867908343738209\\
573.5	0.00785907108188873\\
574	0.00922121973835914\\
574.5	0.00922919486574603\\
575.5	0.00903322638832773\\
576	0.00893971389228913\\
576.5	0.00889620596652905\\
577	0.00876358613621764\\
577.5	0.00867772723670701\\
578	0.00851107381083721\\
578.5	0.00841058356619657\\
579	0.00822691728145977\\
579.5	0.00819850731733904\\
580	0.00801351463659429\\
};
\addlegendentry{relative temperature}

\addplot [color=mycolor2,dashed,line width=\lineWidth]
  table[row sep=crcr]{%
540	0.00440154770378907\\
540.5	0.00419442282742637\\
541	0.00401388010834289\\
541.5	0.00381116891372149\\
542	0.00362649975967567\\
543	0.00326767624574759\\
544.5	0.00272506863332003\\
545	0.00255839509461723\\
545.5	0.00212384901838586\\
546	0.00192840283121467\\
546.5	0.00175811570140424\\
547	0.00157629624191889\\
547.5	0.00139876774416844\\
548	0.00125782045751703\\
548.5	0.00108082503867843\\
549	0.000920438070817091\\
549.5	0.000767648817490632\\
550	0.000632770590766613\\
550.5	0.000495700342000041\\
551	0.000381559134686574\\
551.5	0.000309013837799684\\
552	0.000273543415248378\\
552.5	0.000354350381544407\\
553	0.000706924655411283\\
553.5	0.000848412700535738\\
554	0.000980584300914772\\
554.5	0.00109518110390189\\
555	0.00123788021006589\\
555.5	0.00136434593693466\\
556	0.00150763188548286\\
556.5	0.00164504839386015\\
557	0.00176974347670881\\
557.5	0.00191373772464938\\
558	0.00204215011597655\\
558.5	0.00216045129833323\\
559	0.00230090767827543\\
559.5	0.00241213037427074\\
560	0.00252322227632827\\
560.5	0.0026493859896062\\
561	0.00275267038498641\\
561.5	0.00289661307962941\\
562	0.00300273125905954\\
562.5	0.00174892729658953\\
563	0.00183895971427267\\
563.5	0.00186537350690287\\
564	0.00363270080737022\\
564.5	0.0037504680118282\\
565	0.00386289431701297\\
565.5	0.0020130360969707\\
566	0.00202819545937609\\
566.5	0.00202758373837705\\
567.5	0.00203501785929618\\
568	0.0020266745069579\\
568.5	0.00202152730568086\\
569	0.00202739868768145\\
569.5	0.00200050372706013\\
570	0.0020045877474857\\
570.5	0.00197346672252628\\
571	0.00197670229171115\\
571.5	0.00196059170555475\\
572	0.00194972913059749\\
572.5	0.00192916110048866\\
573	0.00193647566561212\\
573.5	0.00176593526012061\\
574	0.00205763300452025\\
574.5	0.0020630702738804\\
575.5	0.00203014578026186\\
576	0.00201491360644884\\
576.5	0.00201034421781178\\
577	0.00198741374000798\\
577.5	0.00197440324067794\\
578	0.0019450063911205\\
578.5	0.001929561074724\\
579	0.00189735981301338\\
579.5	0.00189703792038601\\
580	0.00186512085119006\\
};
\addlegendentry{supply mass fraction}

\end{axis}
\end{tikzpicture}%